\documentclass[a4paper, fleqn]{article}
\usepackage{sectsty}
\usepackage{amsmath}
\usepackage{amsthm}
\usepackage{dsfont}
\usepackage{pxfonts}

\newtheorem{teo}{Theorem}
\newtheorem{prop}{Proposition}
\newtheorem{coro}{Corollary}
\newtheorem{lem}{Lemma}

\makeatletter
\def\@seccntformat#1{\@ifundefined{#1@cntformat}%
{\csname the#1\endcsname\quad}
{\csname #1@cntformat\endcsname}
}
\def\section@cntformat{\S\,\thesection\quad}
\def\subsection@cntformat{\S\,\thesubsection\quad}
\makeatother

\begin{document}

\def \ni {\noindent}
\allsectionsfont{\mdseries\itshape\centering}

\title{On a ring of modular forms 
related to the Theta gradients map in genus $2$}
\author{A. Fiorentino}
\date{}
\maketitle

\begin{abstract}
\ni The level moduli space $A_g^{4,8}$ is mapped to the projective space by means of gradients of odd Theta functions, such a map turning out no to be injective in the genus $2$ case. In this work a congruence subgroup $\Gamma$ is located between $\Gamma_2(4,8)$ and $\Gamma_2(2,4)$ in such a way the map factors on the related level moduli space $A_{\Gamma}$, the new map being injective on $A_{\Gamma}$. Satake's compactification $\text{Proj}A(\Gamma)$ and the desingularization $\text{Proj}S(\Gamma)$ are also due to be investigated, since the map does not extend to the boundary of the compactification; to aim at this, an algebraic description is provided, by proving a structure theorem both for the ring of modular forms $A(\Gamma)$ and the ideal of cusp forms $S(\Gamma)$.
\end{abstract}

\begin{section}{Basic definitions}
Throughout this paper the symbol $\mathfrak{S}_g$ will stand for the {\it Siegel upper half-plane of degree $g$}, namely the tube domain of complex symmetric $g \times g$ matrices, whose imaginary part is positive definite. A transitive action of the symplectic group $Sp(2g,\mathds{R})$ is defined on  $\mathfrak{S}_g$  by biholomorphic automorphisms:
\begin{equation}
\label{act}
\begin{array}{l} Sp(2g,\mathds{R}) \times \mathfrak{S}_g \longrightarrow \mathfrak{S}_g \\
                   \\
                 \left ( \gamma , \tau \right ) \rightarrow (a_{\gamma} \tau + b_{\gamma}) \cdot (c_{\gamma} \tau + d_{\gamma})^{-1}
\end{array}
\end{equation}
\ni where the generic element of $Sp(2g,\mathds{R})$ is conventionally depicted in a standard block notation as:
\begin{equation*}
 \gamma = \begin{pmatrix} a_{\gamma} & b_{\gamma} \\ c_{\gamma} & d_{\gamma} \end{pmatrix}   \quad \quad \text{with} \, \, a_{\gamma},  b_{\gamma},  c_{\gamma},  d_{\gamma} \,\text{real $g \times g$ matrices} 
\end{equation*}
\ni In particular, $1_{2g}$ and $-1_{2g}$ acts likewise.\\
\ni The action of the {\it Siegel modular group} $\Gamma_g \coloneqq Sp(2g,\mathds {Z}))$ on $\mathfrak{S}_g$ is properly discontinuous; the coset space $A_g \coloneqq {\mathfrak{S}_g}/{\Gamma_g}$ is thus provided with a normal analytic space structure. This space turns out to be remarkably meaningful in the theory of abelian varieties, its points being identified with the classes of isomorphic principally polarized abelian varieties.\\

Subgroups of finite index in $\Gamma_g$ are characterized by containing for some $n \in \mathds{N}$ the level subgroup:
\begin{equation*}
\Gamma_g(n) = \{\, \gamma \in \Gamma_g \mid \gamma \equiv 1_{2g} \, \text{mod} \,  n  \,\}
\end{equation*} which is of finite index itself and normal in $\Gamma_g$; such subgroups are known as {\it congruence subgroups} of the Siegel modular group. Other remarkable families of congruence subgroups this paper will deal with are:

\begin{equation*}
\begin{array}{l}
\Gamma_g(n,2n) \coloneqq \{\, \gamma \in \Gamma_g (n) \mid diag(a_{\gamma}{}^t b_{\gamma} ) \equiv  diag(c_{\gamma}{}^t d _{\gamma} ) \equiv 0 \, \text{mod} \, 2n \,\}\\
\\
\Gamma_g(n,2n,4n) \coloneqq \{\, \gamma \in \Gamma_g (2n,4n) \mid Tr(a_{\gamma}) \equiv g \, \text{mod}\,  n \,\} 
\end{array}
\end{equation*} 

\ni The congruence subgroup $\Gamma_g(n,2n)$ is normal in $\Gamma_g$ whenever $n$ is even; the correspondent level moduli space of principally polarized abelian varieties can be denoted by the symbol $A_g^{n, 2n}$.\\
The quotient $\Gamma_g (2n,4n)/ \Gamma_g(4n,8n)$ is, in fact, supplied with a further peculiar structure (cf. \cite{Igusa3}):

\begin{prop}
\label{generalgammabasis}
$\Gamma_g(2n,4n) / \Gamma_g(4n,8n)$ is a $g(2g+1)$-dimensional vector space on $\mathds{Z}_2$; in particular,  a basis for $\Gamma_g(2,4) / \Gamma_g(4,8)$ is given by the representatives of the following elements:
\begin{equation*}
\begin{array}{lll}
A_{ij} \coloneqq \begin{pmatrix} a_{ij} & 0 \\ 0 & {}^t a_{ij}^{-1}\end{pmatrix} & \quad (1 \leq i,j < g) & \quad \quad \quad a_{ij} \coloneqq \left \{ \begin{array}{ll}1_g+2\tilde{O}_{ij} & \text{if} \, \,  i \neq j \\ 1_g - 2\tilde{O}_{ij} & \text{if} \, \,  i = j \end{array} \right. \\
& \\
B_{ij} \coloneqq \begin{pmatrix} 1_g & 2\tilde{O}_{ij} + 2\tilde{O}_{ji} \\ 0 & 1_g \end{pmatrix} & \quad (1 \leq i < j \leq g) & \\
& \\
B^2_{ii} \coloneqq \begin{pmatrix} 1_g & 2\tilde{O}_{ii}  \\ 0 & 1_g \end{pmatrix}^2 & \quad (1 \leq i \leq g) & \\
& \\
C_{ij} \coloneqq {}^t B_{ij} & \quad (1 \leq i <  j \leq g) & \\
\\
C^2_{ii} \coloneqq {}^t B^2_{ii} & \quad (1 \leq i \leq g) &\\
\\
-1_{2g}
\end{array}
\end{equation*}\\
where $\tilde{O}_{ij}$ stand for the $g \times g$ matrix, whose $(h,k)$ coordinates are $\tilde{O}_{ij}^{(hk)}=\delta_{ih}\delta_{jk}$. 
\end{prop}

\ni The quotient $\Gamma_2(2,4)/ \{ \pm \Gamma_2(4,8) \}$ this work is concerned with will be conventionally denoted by the symbol $G$, the group of characters of $G$ being meant by $\hat{G}$.\\

A {\it Siegel modular form} of weight $k \in \mathds{Z}^+$ with respect to a congruence subgroup $\Gamma \subset \Gamma_g$ is a holomorphic \footnote{When $g=1\,$  $f$ must be also holomorphic on $\infty$, this condition being turned into a redundant one by (\ref{modularity}) when $g > 1$.} function $f: \mathfrak{S}_g \rightarrow \mathds{C}$, satisfying the so-called modularity condition:
\begin{equation}
\label{modularity}
f(\gamma \tau) = \det{(c_{\gamma} \tau + d_{\gamma})}^k f(\tau) \quad \quad \quad \forall \gamma \in \Gamma , \quad \forall \tau \in \mathfrak{S}_g
\end{equation}
\ni Henceforward Siegel modular forms will be simply referred to as modular forms; the ring $A(\Gamma)$ of modular forms with respect to $\Gamma$ is positively graded by the weights, namely $A(\Gamma) = \bigoplus_{k \geq 0}A_k(\Gamma)$. \\
For $1 \leq k \leq g$ the {\it Siegel operator} $\Phi_{g, k}:  A(\Gamma_g) \rightarrow A(\Gamma_k)$ acting on modular forms is well defined by setting:
\begin{equation*}
\label{Siegeloperatorinformal}
\Phi_{g, k}(f) (\tau) = \lim_{\lambda \rightarrow \infty} f\begin{pmatrix} \tau & 0 \\ 0 & i \lambda  \end{pmatrix} \quad \quad \forall f \in A(\Gamma_g) , \, \, \, \forall \tau \in \mathfrak{S}_{k}
\end{equation*} and is seen to preserve the weight. By means of such an operator one is allowed to describe the ideal of the so-called {\it cusp forms}, namely the modular forms that vanish on the boundary of Satake's compactification:
\begin{equation*}
\begin{array}{ccc}
S(\Gamma) \coloneqq \bigoplus_{k \geq 0}S_k(\Gamma) & \text{where} &
S_k(\Gamma) \coloneqq \{ \, f \in A_k(\Gamma) \, | \, \Phi_{g, g-1}(\gamma^{-1}|_kf)=0 \quad \forall \gamma \in \Gamma_g \, \}
\end{array}
\end{equation*}\\

{\it Riemann Theta functions with characteristics} play an essential role in the construction of several modular forms; for each $m=(m' , m'') \in \mathds{Z}^g \times \mathds{Z}^g$, they are defined as holomorphic functions on $\mathfrak{S}_g \times \mathds{C}^g$ by the series:
\begin{equation*}
\theta_{m}(\tau , z) \coloneqq \sum_{n \in \mathds{Z}^g} \text{exp} \left \{ {}^t \left ( n + \frac{m'}{2} \right ) \tau \left ( n + \frac{m'}{2} \right ) + 2 {}^t \left ( n + \frac{m'}{2} \right ) \left ( z + \frac{m''}{2} \right ) \right \}
\end{equation*} where $\text{exp}(z)$ stands for the function $e^{\pi iz}$. Since $\theta_{m+2n}(\tau , z) = (-1)^{{}^tm'n''} \theta_{m} (\tau , z )$, Theta functions can be reduced to the only ones related to {\it $g$-characteristics}, namely vector columns $\begin{bmatrix} m' \\ m'' \end{bmatrix}$ with $m', m'' \in \mathds{Z}_{2}^g$. The set of $g$-characteristics will be denoted by the symbol $C^{(g)}$; needless to say, each $g$-characteristic $m$ satisfies $m + m =0$. For each $m \in C^{(g)}$ the function $\theta_m: \mathfrak{S}_g \rightarrow \mathds{C}$, defined by $\theta_m (\tau) \coloneqq \theta_m(\tau , 0)$, is known as the {\it Theta constant} with $g$-characteristic $m$ (or simply with characteristic $m$, when no misunderstanding is allowed). By introducing the parity function  $e(m)=(-1)^{^t m' m''}$, characteristics can be classified into {\it even} and {\it odd} ones, respectively if $e(m)=1$ or $e(m)=-1$; Theta constants related to even characteristics are thus plainly checked to be the only non-vanishing ones. Since there are $2^{g-1}(2^g+1)$ even $g$-characteristics and $2^{g-1}(2^g-1)$ odd $g$-characteristics, there exist exactly $2^{g-1}(2^g+1)$ non-vanishing Theta constants for each $g \geq 1$ . \\
\ni To review how Theta constants transforms under $\Gamma_g$, an action on the set $C^{(g)}$ has to be introduced by setting:
\begin{equation}
\label{action}
\gamma \begin{bmatrix} m' \\ m'' \end{bmatrix} \coloneqq \left [ \begin{pmatrix} d_{\gamma} &  -c_{\gamma} \\ - b_{\gamma} &  a_{\gamma} \end{pmatrix} \begin{pmatrix} m' \\ m'' \end{pmatrix} + \begin{pmatrix} diag(c_{\gamma}{}^t  d_{\gamma} ) \\ diag( a_{\gamma}{}^t  b_{\gamma})  \end{pmatrix} \right ] \text{mod}\, 2 
\end{equation}This action is easily seen to preserve the parity, $C^{(g)}$ being decomposed into even and odd characteristics by it.\\
\ni With reference to the actions in (\ref{act}) and (\ref{action}), a transformation law for Riemann Theta functions holds:
\begin{equation}
\label{generaltranformtheta}
\begin{array}{l}
\theta_{\gamma m} (\gamma \tau \, , {}^t(c_{\gamma} \tau + d_{\gamma})^{-1} z) = \Phi(m,\gamma, \tau, z) {det(c_{\gamma} \tau + d_{\gamma})}^{\frac{1}{2}} \theta_m (\tau, z) \\
\\
\forall \gamma \in \Gamma_g, \quad \forall m \in C^{(g)},  \quad \forall \tau \in \mathfrak{S}_g, \quad \forall z \in \mathds{C}^g
\end{array}
\end{equation}
where $\det{(c_{\gamma} \tau + d_{\gamma})}^{\frac{1}{2}}$ stands for the branch of the rooth whose sign turns to be positive when $Re\tau = 0$, and the function $\Phi$ can be split into a suitable product of two factors, one of them depending on the sole variable $\gamma$:
\begin{equation*}
\Phi(m,\gamma, \tau, z) = \kappa(\gamma) \exp \left \{ \frac{1}{2} {}^t z \,  [(c_{\gamma} \tau + d_{\gamma})^{-1} c_{\gamma}] \, z + 2\phi_m(\gamma) \right \}
\end{equation*}
\ni with
\begin{equation*}
\label{generalphi}
\begin{split}
\phi_m(\gamma)  = & -\frac{1}{8} ({}^t m ' {}^tb_{\gamma} d_{\gamma} m' + {}^t m '' {}^t a_{\gamma} c_{\gamma} m'' - 2{}^t m ' {}^t b_{\gamma} c_{\gamma} m'') + \\
                                         &  -\frac{1}{4} {}^tdiag(a_{\gamma}{}^t b_{\gamma})(d_{\gamma} m ' - c_{\gamma} m '')
\end{split}
\end{equation*}

\ni An outstanding peculiarity of the function $\kappa$ is that $\kappa^2$ is a character of $\Gamma_g(1,2)$. On $\Gamma_g(2)$ the function $\kappa^2$ admits, in particular, a simple explicit expression (cf. \cite{Igusa3}, \cite{Igusa10} and \cite{CuspForms}):
\begin{equation}
\label{k}
\begin{array}{l}
{\kappa (\gamma)}^{2} = e ^{ \frac{\pi}{2} Tr(a_{\gamma}-1_g) i} \quad \quad \forall \gamma \in \Gamma_g(2)
\end{array}
\end{equation}

\ni $\kappa^2$ is, therefore, a character of $\Gamma_g(2,4)/\{\pm \Gamma(4,8) \}$ when $g$ is even, the function beeing trivial on $\Gamma(4,8)$ and well defined on $\Gamma_g(2)/\{\pm 1_{2g} \}$ whenever $g$ is even.

\ni By setting: 
\begin{equation*}
\label{chiexpression}
\chi_m (\gamma) \coloneqq  \Phi (m, \gamma, \tau, 0) = e^{2 \pi i \phi_m(\gamma)}
\end{equation*} the law $\theta_{\gamma m} (\gamma \tau) = \kappa(\gamma) \chi_m(\gamma) {det(c_{\gamma} \tau + d_{\gamma})}^{\frac{1}{2}} \theta_m (\tau)$ is gained, (\ref{generaltranformtheta}) being applied to Theta constants; when $\gamma \in \Gamma_g(2)$ this yields the transformation law:
\begin{equation}
\label{transformcostanthetagamma2n4n}
\theta_{m} (\gamma \tau) = \kappa(\gamma) \chi_m(\gamma) \det{(c_{\gamma} \tau + d_{\gamma})}^{\frac{1}{2}} \theta_m (\tau)
\end{equation}
\ni for $\Gamma_g(2)$ acts trivially on $C^{(g)}$.\\

\ni As concerns $\chi_m$, one has (cf. \cite{Igusa3} or \cite{SM}):

\begin{lem}
\label{chimonbasis}
$\chi_m$ is a character of $\Gamma_g(2,4)/\Gamma_g(4,8)$ \footnote{$\chi_m^2$ is in fact a character of $\Gamma_g(2)/\Gamma_g(2,4)$ (cf. \cite{SM}). } and, with reference to the basis described in Proposition \ref{generalgammabasis}, one has:\\
\begin{equation*}
\begin{array}{ccc}
\begin{array}{l}
\chi_m(B^2_{ii})=(-1)^{m_i'^2} \\
\\
\chi_m(B_{ij})=(-1)^{m_i'm_j'}\\
\\
\chi_m(C^2_{ii})=(-1)^{m_i''^2}\\
\\
\chi_m(C_{ij})=(-1)^{m_i'''m_j''}
\end{array} & & 
\begin{array}{l}
\chi_m(A_{ij})=(-1)^{m_i'm_j''} \\
\\
\chi_m(-1_{2g})= \left \{ \begin{array}{ll} 1 & \text{if} \, \, $m$ \, \, \text{is even} \\
                                              &                                    \\
                                           -1 & \text{if} \, \, $m$ \, \, \text{is odd}
\end{array} \right.
\end{array}
\end{array}
\end{equation*} where $m'_i$ and $m''_i$ denote respectively the i-th coordinate of $m'$ and the i-th coordinate of $m''$ in $m= \begin{bmatrix} m' \\ m'' \end{bmatrix}$ .
\end{lem}

\ni Since $\chi_m$ is trivial on $\Gamma_g(4,8)$, (\ref{transformcostanthetagamma2n4n}) implies that the product $\theta_m \theta_n$ of two Theta constants is a modular form of weight $1$ with respect to $\Gamma_g(4,8)$.\\

\ni A useful criterion of modularity with respect to $\Gamma_g(2,4)$ can be more generally stated for products of even sequences of Theta constants (cf. \cite{SM}):
\begin{prop}
\label{modular}
Let $M=(m_1, \dots , m_{2k})$ be a sequence of even characteristics. The product ${\theta}_{m_1} \cdots {\theta}_{m_{2k}}$ is a modular form with respect to $\Gamma_g(2,4)$ if and only if:
\begin{equation}
\label{modularcriterium}
M\,{}^tM \equiv k\begin{pmatrix} 0 & 1 \\ 1 & 0 \end{pmatrix} \, \text{mod}  \, 2
\end{equation}
\end{prop}

\end{section}

\begin{section}{Further notions about $2$-characteristics}
\label{mainnotations}
To prove the results this work aims at, some combinatorial features of $2$-characteristics are demanded; this section is thus intended to point out these peculiarities.\\
\ni A special notation for the case $g=2$ will be introduced here for the sake of simplicity and referred to throughout the paper. The symbol $C_k$ will denote the parts of the set of even $2$-characteristics whose cardinality is $k$, while $\tilde{C}_k$ will mean the parts of the set of odd $2$-characteristics whose cardinality is $k$; $C_1$ and  $\tilde{C}_1$, in particular, will stand respectively for the set of even $2$-characteristics and for the set of odd $2$-characteristics. The cardinality of a set $C$ will be generally denoted in this paper by means of the symbol $|C|$.\\

The group $\Gamma_2/\Gamma_2(2) \cong S_6$ acts both on $C_1$ and $\tilde{C}_1$, the action of $\Gamma_2(2)$ being trivial; by focusing on the action naturally defined on the sets $C_k$ for each $k$, the following decompositions into orbits arises (cf. \cite{Structure}):

\begin{enumerate}
\item $C_2$ is a unique orbit.
\item $C_3= C_{3}^{-} \cup C_{3}^{+}$ ($|C_{3}^{-}| = |C_{3}^{+}|=60$), where:
\begin{equation*}
\begin{array}{l}
 C_{3}^{-}= \{ \, \{m_1, m_2, m_3\} \in C_3 \, \, | \, \, m_1 + m_2 + m_3 \in \tilde{C}_1 \} \\
\\
 C_{3}^{+}= \{ \, \{m_1, m_2, m_3\} \in C_3 \, \, | \, \, m_1 + m_2 + m_3 \in C_1 \}
\end{array}
\end{equation*}

\item $C_4 = C_4^+ \cup C_4^* \cup C_4^-$  ($|C_{4}^{-}| = |C_{4}^{+}|=15$ , $|C_4^*| = 180$), where:
\begin{equation*}
\begin{array}{l}
 C_{4}^{-}= \{ \, \{m_1, \dots, m_4\} \in C_4 \, \, | \, \, \{m_i, m_j, m_k\} \in C_{3}^{-}\, \, \, \forall   \{m_i, m_j, m_k\} \subset  \{m_1, \dots, m_4\} \, \} \\
\\
 C_{4}^{+}= \{ \, \{m_1, \dots , m_4\} \in C_4 \, \, | \, \, \{m_i, m_j, m_k\} \in C_{3}^{+} \, \, \, \forall   \{m_i, m_j, m_k\} \subset  \{m_1, \dots, m_4\} \, \}
\end{array}
\end{equation*}

\item $C_5 = C_5^+ \cup C_5^* \cup C_5^-$  ($|C_{5}^{-}| = |C_{5}^{+}|=90$ , $|C_5^*| = 72$), where:
\begin{equation*}
\begin{array}{l}
C_{5}^{-} = \{ \, \{m_1, \dots m_5\} \in C_5 \, \, | \, \, \{m_1, \dots m_5\} \, \text{contains a unique element of} \, C_{4}^{-} \} \\
\\
C_{5}^{+} = \{ \, \{m_1, \dots m_5\} \in C_5 \, \, | \, \, \{m_1, \dots m_5\} \, \text{contains a unique element of} \, C_{4}^{+} \} \}
\end{array}
\end{equation*}

\item For $k > 5$ the decomposition is likewise given by complementary sets:
\begin{equation*}
\begin{array}{l}
 C_{k}^{-}= \{ \, \{m_1, \dots, m_k\} \in C_k \, \, | \, \,\{m_1, \dots, m_k\}^c \in {C_{10-k}}^{+} \} \\
\\
 C_{k}^{+}= \{ \, \{m_1, \dots, m_k\} \in C_k \, \, | \, \,\{m_1, \dots, m_k\}^c \in {C_{10-k}}^{-} \}
\end{array}
\end{equation*}
\end{enumerate}

\ni A complete description of these orbits is displayed in a diagram in \cite{Structure}; here, some needed features will be briefly reviewed.
\begin{lem}
\label{l1}
For each $m_1, m_2 \in  C_1$ distinct, four characteristics $n_1, n_2, n_3, n_4$ in the complementary set $\{m_1, m_2\}^c$ are such that $\{m_1, m_2,  n_i \} \in C_3^{-}$ and the other four $h_1, h_2, h_3, h_4$ such that $\{m_1, m_2,  h_i \} \in C_3^{+}$.
\end{lem} 

\begin{lem}
\label{l0}
Let $\{ m_1, m_2,  m_3 \} \in C_3^{-}$. There is exactly one characteristic $n \in C_1$ such that $\{m_1, m_2,  m_3,  n \} \in C_4^{-}$.
\end{lem} 

\begin{lem}
\label{l0+}
Let $\{ m_1, m_2,  m_3 \} \in C_3^{+}$. There is exactly one characteristic $n \in C_1$ such that $\{m_1, m_2,  m_3,  n \} \in C_4^{+}$, namely $n = m_1 + m_2 + m_3$.
\end{lem} 

\begin{coro}
\label{c4+}
$C_4^{+}=\{\,\{m_1, m_2, m_3, m_4\} \in C_4 \mid \sum_{i=1}^4 m_i = 0\,\}$
\end{coro}

\begin{coro}
\label{c1}
Let $h_1, h_2 \in  C_1$ be distinct. There are exactly two elements in $C_4^{-}$ containing $\{h_1, h_2 \}$. 
\end{coro}

\ni More precisely, one has:

\begin{lem}
\label{l1t}
If $\{m_1,m_2,h,k \}, \{m_3,m_4,h,k \} \in C_4^{-}$, then $\{m_1,m_2,m_3,m_4 \} \in C_4^{-}$.
\end{lem}

\begin{lem}
\label{l2t}
Let $\{m_1,m_2,m_3,n \}, \{m_4,m_5,m_6,n \} \in C_4^{-}$ be such that $m_i, n$ are all distinct;  then $\{h_1, h_2, h_3, n \} \in C_4^{-}$, where $\{h_1, h_2, h_3\} = \{m_1,m_2,m_3, m_4,m_5,m_6, n \}^c$
\end{lem}

\ni Moreover:
\begin{prop}
\label{c5-c4-}
If $M \notin C_{5}^{-}$, $M$ does not contain any element of $C_4^{-}$.
\end{prop}

\begin{prop}
\label{c6-c5-}
$M \in C_6^{-}$ if and only if $M$ contains exactly six elements of $ C_{5}^{-}$.
\end{prop}

\begin{coro}
\label{c6-sum0}
$M \in C_6^{-}$ if and only if $\sum_{m \in M} m = 0$. 
\end{coro}

\begin{coro}
\label{l00c}
$M \in C_6^{-}$ (respectively $M \in C_6^{+}$) if and only if $M$ does not contain any element of $C_4^{+}$ ($C_4^{-}$).
\end{coro}

\ni Elements $M_1, M_2$ belonging to $C_4^{+}$ or $C_6^{-}$ clearly satisfies:

\begin{equation}
\label{c6symc6}
0 = \sum_{m \in M_1} m + \sum_{m \in M_2} m = \sum_{m \in M_1 \vartriangle M_2} m 
\end{equation}

\ni where the symbol $\vartriangle$ stands for the classical symmetric difference between sets, namely $M_i \vartriangle M_j \coloneqq (M_i \cup M_j) - (M_i \cap M_j)$. The following Propositions describe the behaviour of the symmetric difference between elements belonging to such distinguishing orbits.

\begin{prop}
\label{pA}
Let $M_1,M_2 \in C_4^{-}$. Then $M_1 \cap M_2  \neq \emptyset$ and:
\begin{equation*}
\begin{array}{llll}
1) & M_1 \vartriangle M_2 \in C_6^{+}& if & |M_1 \cap M_2| = 1 \\
2) & M_1 \vartriangle M_2 \in C_4^{-} & if & |M_1 \cap M_2| = 2 \\
3) & M_1=M_2                                          & if & |M_1 \cap M_2| > 2 \\
\end{array} 
\end{equation*}
\end{prop}

\begin{proof}
Lemma \ref{l00c} implies $M_2 \nsubset M_1^c$, hence $M_1 \cap M_2  \neq \emptyset$. Then the cases in the statement are the only possible ones.\\
\ni $1)$ If $M_1 = \{m_1,m_2,m_3,n\}$ and $M_2 = \{m_4,m_5,m_6,n\}$ with $m_i , n$ all distinct, then $\, M_1 \vartriangle M_2=   \{m_1,m_2,m_3,m_4,m_5,m_6\} \in C_6^{+}$ by Lemma \ref{l2t}.\\
\ni $2)$ If $M_1 = \{m_1,m_2,h,k\}$ and $M_2 = \{m_3,m_4,h,k\}$, $\, M_1 \vartriangle M_2=   \{m_1,m_2,m_3,m_4\} \in C_4^{-}$ by Lemma \ref{l1t}.\\
\ni $3)$ It is due to Lemma \ref{l0}.
\end{proof}

\begin{prop}
\label{pB}
Let $M_1,M_2 \in C_6^{+}$. Then $|M_1 \cap M_2| \geq 3$ and:
\begin{equation*}
\begin{array}{llll}
1) & M_1 \vartriangle M_2 \in C_6^{+} & if & |M_1 \cap M_2| = 3 \\
2) & M_1 \vartriangle M_2 \in C_4^{-}& if &  |M_1 \cap M_2| = 4\\
3) & M_1 = M_2                                        & if & |M_1 \cap  M_2| >4 
\end{array} 
\end{equation*}
\end{prop}

\begin{proof}
$|M_1 \cap M_2| \geq 2$ trivially because $|C_1|=10$; moreover:
\begin{equation}
\label{m1m26partition}
|M_1 \cap M_2| - |M_1^c  \cap M_2^c | = 2
\end{equation}
\ni since $12=|M_1|+|M_2|$. Hence, $|M_1 \cap M_2| > 2$, for $M_1^c \cap M_2^c \neq \emptyset$ by Proposition \ref{pA}; the cases described in the statement are therefore the only possible ones. \\
$1)$ If $M_1 = \{m_1,m_2,m_3,h,k,l\}$ and $M_2 = \{m_4,m_5,m_6,h,k,l\}$, then one has $M_1^c = \{m_4,m_5,m_6,n\} \in C_4^{-}$ and $M_2^c = \{m_1,m_2,m_3,n\} \in C_4^{-}$, where $n \neq m_i,h,k,l$; Hence, $M_1 \vartriangle M_2= \{h,k,l,n\}^c \in C_6^{+}$ by Lemma \ref{l2t}.\\
$2)$ If $M_1 = \{m_1,m_2,n,h,k,l\}$ and $M_2 = \{m_3,m_4,n,h,k,l\}$, then one has $M_1^c = \{m_3,m_4, i, j\} \in C_4^{-}$ and $M_2^c = \{m_1,m_2,i,j\} \in C_4^{-}$, where $i,j \neq m_i,n,h,k,l$. Hence, $M_1 \vartriangle M_2=   \{m_1,m_2,m_3,m_4\} \in C_4^{-}$ by Lemma \ref{l1t}. \\
$3)$ If $|M_1 \cap M_2| >4$, then $|M_1^c \cap M_2^c| > 2$ by (\ref{m1m26partition}). Hence, $M_1^c=M_2^c$ by Proposition \ref{pA}, and, therefore, $M_1=M_2$.
\end{proof}

\begin{prop}
\label{pC}
Let $M_1 \in C_6^{+}$ and $M_2 \in C_4^{-}$.  If $M_1^c \neq M_2$ the only possible cases are:
\begin{equation*}
\begin{array}{llll}
1) & M_1 \vartriangle M_2 \in C_4^{-}& if & |M_1 \cap M_2| = 3 \\
2) & M_1 \vartriangle M_2 \in C_6^{+} & if & |M_1 \cap M_2| = 2
\end{array} 
\end{equation*}
\end{prop}

\begin{proof}
Obviously $|M_1 \cap M_2| \leq 4$; moreover,  Lemma \ref{l00c} implies $M_2 \nsubset M_1$, hence $|M_1 \cap M_2| \leq 3$. If $M_1^c \neq M_2$, Lemma \ref{l0} implies $|M_1^c \cap M_2|<3$, hence $|M_1 \cap M_2| >1$. The only possible cases when $M_1^c \neq M_2$ are, therefore, the ones described in the statement.\\
$1)$ If $M_1=\{m_1,m_2,m_3,h,k,l\}$ and $M_2 = \{m_4,h,k,l\}$, $M_2^c = \{m_1,m_2,m_3,r,s,t\} \in C_6^{+}$ with $r,s,t \neq m_i,h,k,l$; hence, by Proposition \ref{pB}:
\begin{equation*}
M_1 \vartriangle M_2=   \{m_1,m_2,m_3,m_4\} = \{h,k,l,r,s,t\}^c=(M_1 \vartriangle M_2^c)^c \in C_4^{-}
\end{equation*}
$2)$ If $M_1=\{m_1,m_2,m_3,m_4,h,k\}$ and $M_2 = \{m_5,m_6,h,k\}$, $M_1^c = \{m_5,m_6,i,j\} \in C_4^{-}$ with $i,j \neq m_i,h,k$; hence, by Proposition \ref{pA}:
\begin{equation*}
M_1 \vartriangle M_2=   \{m_1,m_2,m_3,m_4,m_5,m_6\} = \{i,j,h,k\}^c=(M_1^c \vartriangle M_2)^c \in C_6^{+}
\end{equation*}

\end{proof}

\end{section}

\begin{section}{The Theta gradients map}

As concerns the gradients:
\begin{equation*}
\text{grad}^0_{z} \theta_n \coloneqq \text{grad}_{z} \theta_n {\mid}_{z=0} = \left ( \frac {\partial}{\partial z_1 } \theta_n {\mid}_{z=0} , \dots , \frac {\partial}{\partial z_g } \theta_n {\mid}_{z=0}  \right )
\end{equation*} the only non trivial ones are those which are related to odd characteristics; a peculiar transformation law holds for them:
\begin{equation}
\label{transformgradient}
\begin{array}{l}
\text{grad}^0_{z} \theta_n (\gamma \tau) = \det{(c_{\gamma} \tau + d_{\gamma})}^{\frac{1}{2}} (c_{\gamma} \tau + d_{\gamma}) \cdot \text{grad}^0_{z} \theta_n (\tau) \\
\\
\forall \gamma \in \Gamma_g(4,8), \quad \forall \tau \in \mathfrak{S}_g
\end{array}
\end{equation}
\ni meaning that gradients of odd Theta functions can be regarded as modular forms with respect to $\Gamma_g(4,8)$ under the representation $T_0(A) \coloneqq {\det{(A)}}^{1/2} A$ \footnote{The classical modular forms of weight $k$, which are the only ones this work is concerning with, are indeed modular forms under the representation $T(A) \coloneqq {\det{(A)}}^{k}$}.
\ni One is thus allowed by (\ref{transformgradient}) to define a map on the space  $A_g^{4,8}$:
\begin{equation*}
\begin{array}{c}
\mathds{P}grTh: A_g^{4,8} \longrightarrow {\overbrace{{\mathds{C}}^g \times \cdots \times {\mathds{C}}^g}^{2^{g-1}(2^g-1) \text{times}}}/T_0(Gl(g,\mathds{C})) \\
                                                   \\
                  \quad \quad \quad \tau \longrightarrow \left \{ {\mbox{grad}_{z} \theta_n} {\mid}_{z=0}  \right \}_{n \, \text{odd}}\end{array}
\end{equation*} whose range lies in the Grassmannian $Gr_{\mathds{C}}(g,2^{g-1}(2^g-1))$ by Lefschetz's theorem (cf. \cite{SalvatiManniPhd}), the Jacobian determinants being the  Pl\"{u}cker coordinates of this map. \\

The investigation around the injectivity of this map is strictly related to the problem of recovering plane curves form their tangent hyperplanes. In fact, if $\phi: C \rightarrow \mathds{P}^{g-1}$ denotes the canonical map of a smooth curve $C$ of genus $g$, and $\tau_C = [J(C)] \in J_g \subset A_g$ the correspondent point in the locus of Jacobians $J_g$, a hyperplane $H \subset \mathds{P}^{g-1}$ being tangent to the canonical curve $\phi(C)$ in $g-1$ points cut a divisor on $\phi(C)$, which is the zero locus of one of the $2^{g-1}(2^g - 1)$ Riemann Theta functions with odd characteristics $\theta_n (z) = \theta_n(\tau_C, z)$; on the converse, each Riemann Theta function with odd characteristic related to the curve $C$ determines such a hyperplane $H \subset \mathds{P}^{g-1}$, whose direction is given by the gradient of the correspondent Riemann Theta function in zero. Hence, $\mathds{P}grTh$ maps an element $\tau \in J_g^{4,8} \subset A_g^{4,8}$ to an ordered set of all the hyperplanes tangent in $g-1$ points; it is, therefore, related to the map sending $\tau \in J_g$ to the set of all the hyperplanes tangent in $g-1$ points, which has been proved to be injective in \cite{capsern1} and \cite{capsern2}; such a map factors through $\mathds{P}grTh$ where both are defined (cf. \cite{mainarticle}). \\
In \cite{mainarticle} the map $\mathds{P}grTh$ is proved to be generically injective on $A_g^{4,8}$ when $g \geq 3$ and injective on tangent spaces when $g \geq 2$; it was also conjectured to be injective whenever $g \geq 3$, albeit this has not been proved yet. \\
\ni Regarding he case $g=2$, the $6$ odd characteristics yield $15$ distinct Jacobian determinants. By setting:
\begin{equation*}
D(n_i, n_j) (\tau) \coloneqq \frac{1}{\pi^2} \begin{vmatrix} \, \, \frac {\partial}{\partial z_1 } \theta_{n^{(i)}}|_{ \, z=0} \, (\tau)  & \frac {\partial}{\partial z_2 } \theta_{n_j}|_{\, z=0} \, (\tau) \\
 & \\ 
\frac {\partial}{\partial z_1 } \theta_{n_i}|_{\, z=0} \, (\tau) & \frac {\partial}{\partial z_2 } \theta_{n_j}|_{\, z=0} \, (\tau)  \, \,  \end{vmatrix}
\end{equation*} the following transformation law arises from the modular properties:
\begin{equation}
\label{trasfDM2}
\begin{array}{l}
D(N)(\gamma \tau) =\kappa (\gamma)^{2} \chi_{N}(\sigma) {det(c_{\gamma} \tau + d_{\gamma})}^{2} D(N)(\tau)\\
\\
\forall \, \tau \in \mathfrak{S}_2 \quad \quad \forall \, \gamma \in \Gamma_2(2,4) \quad \quad  \forall \, N=\{n_1, n_2 \} \in \tilde{C}_2 
\end{array}
\end{equation} where $\chi_N = \chi_{n_1}\chi_{n_2}$. Since $\kappa (\gamma)^{2} \chi_{N}(\gamma)$ is a sign and  $\Gamma_2(2,4) / \Gamma_2(4,8)$ contain $2^{10}$ elements by Proposition \ref{generalgammabasis}, there are $\sum_{k=0}^6\begin{pmatrix} 6 \\ k\end{pmatrix}=64$ possible image points for each set $\{\gamma \tau_0\}_{[\gamma] \in \Gamma_2(2,4) / \Gamma_2(4,8)}$ through the map:

\begin{equation*}
\begin{array}{l}
                           \mathds{P}grTh_2:A^{4,8} \longrightarrow \mathds{P}^{14} \\
                                                      \\
                      \quad \quad \quad \tau \longrightarrow [ \, D(N_1)(\tau), \cdots D(N_{15}) (\tau) \, ]
\end{array}
\end{equation*} 
\ni Hence,  $\mathds{P}grTh_2$ can not be injective, although it is finite of degree $16$ (cf. \cite{mainarticle}). However, a suitable congruence subgroup $\Gamma$ is located between $\Gamma_2(4,8)$ and $\Gamma_2(2,4)$ in such a way that $\mathds{P}grTh_2$ is still well defined on the correspondent level moduli space $A_{\Gamma} \coloneqq \mathfrak{S}_2/\Gamma$, being also injective on it. The next section will be intended to describe such a remarkable group.\\

\ni This section concludes by noting that (\ref{trasfDM2}) implies each $D(N)$ is a modular form of weight $2$ with respect to $\Gamma_2(4,8)$. It also yields a criterion for products of $D(N)$, which is similar to the one described in Proposition \ref{modular}: $D(N_1) \cdots D(N_h)$ is a modular form with respect to $\Gamma_2(2,4)$$\,$\footnote{The criterion still holds, indeed, in the general case, the $D(N)$ being modular forms of weight $\frac{1}{2}(g+2)$.} if and only if:
\begin{equation}
\label{Nmodular}
N\,{}^tN \equiv h\begin{pmatrix} 0 & 1 \\ 1 & 0 \end{pmatrix} \, \, \text{mod} \, 2
\end{equation} with $N=(N_1, \dots n_h)$.

\end{section}

\begin{section}{The congruence subgroup $\Gamma$}
\label{oddchar}
An {\it a priori} description for $\Gamma$ is indeed provided by (\ref{trasfDM2}); denoting by $\chi_{N_1}, \dots \chi_{N_{15}}$ the characters appearing in the transformations of  the fifteen Jacobian determinants under $\Gamma_2(2,4)$, it must be:

\begin{equation}
\label{condicio}
\Gamma \coloneqq \Gamma^{(1)} \cup \Gamma^{(-1)} 
\end{equation}

\ni where:

\begin{equation*}
\begin{array}{lll}
\Gamma^{(1)} & \coloneqq & \{\,\gamma \in \Gamma_2(2,4) \, \mid \kappa(\gamma)^{2} \chi_{N_i}(\gamma) = 1 \quad \forall i= 1, \dots , 15 \,\}\\
& & \\
\Gamma^{(-1)} &  \coloneqq & \{\,\gamma \in \Gamma_2(2,4) \, \mid \kappa(\gamma)^{2} \chi_{N_i}(\gamma) = -1 \quad \forall i= 1, \dots , 15 \,\}
\end{array}
\end{equation*} 

\ni The set $\Gamma$, as defined in (\ref{condicio}), is clearly a subgroup of the Siegel modular group $\Gamma_2$, due to the properties of $\kappa^2$ and $\chi_{N_i}$; moreover, $\Gamma_2(4,8) \subset \Gamma$, all these characters being trivial on $\Gamma_2(4,8)$. The expression in (\ref{condicio}) defines, therefore, a congruence subgroup fulfilling all the requirements.\\
\ni The next step is to refine the provisional definition of $\Gamma$, by detecting which elements of $\Gamma_2(2,4)$ actually belong to it; this is the purpose the remaining part of this section is mostly  concerned with.

\begin{prop}
\label{propcondicio}
\begin{equation}
\label{gamma}
\bigcap_{i= 1}^{15} \text{Ker}\chi_{N_i} = \Gamma = \Gamma^{(1)}
\end{equation}
\end{prop}

\begin{proof}
\ni Since $\bigcap_{i= 1}^{15} \text{Ker}\chi_{N_i} \subset \Gamma$, the sole reverse inclusion has to be shown to prove the first identity. Let, thus, $\gamma \in \Gamma$; if $\chi_{N_i}(\gamma)=-1$  for each $i=1, \dots , 15$, an absurd statement clearly turns up:
\begin{equation*}
-1 = \chi_{(n, n_i)}(\gamma) \, \chi_{(n, n_j)} (\gamma) \, \chi_{(n, n_k)}(\gamma)= \chi_{(n_i, n_j)}(\gamma) \, \chi_{(n, n_k)}\,(\gamma)=1
\end{equation*} Therefore, by (\ref{condicio}) $\chi_{N_i}(\gamma)=1$ for each $i=1, \dots , 15$; hence, $\gamma \in \bigcap_{i= 1}^{15} \text{Ker}\chi_{N_i}$.\\
\ni To prove the second identity in the statement, one has to show that $k^2(\gamma)=1$ whenever $\gamma \in \Gamma$.  By applying the criterion (\ref{Nmodular}), the products:
\begin{equation*}
D \coloneqq D(n_1, n_2)D(n_3, n_4)D(n_5, n_6) \quad \quad \quad n_1, \dots, n_6 \quad \text{all distincts}
\end{equation*}
\ni are checked to be modular forms with respect to $\Gamma_2(2,4)$. However, by (\ref{trasfDM2}):
\begin{equation*}
D(\gamma \tau) =k^2(\gamma) \chi_{n_1}\cdots \chi_{n_6}det(c_{\gamma}\tau+d_{\gamma})^6D(\tau)
\end{equation*}
\ni and, since each $\chi^2_n$ is trivial on $\Gamma_2(2,4)$, one has therefore:
\begin{equation*} 
k^2(\gamma)=\prod_{i=1}^6\chi_{n_i}(\gamma)=\chi_{(n_1,n_2)}(\gamma)\, \chi_{(n_3,n_4)}(\gamma)\, \chi_{(n_5,n_6)}(\gamma)  \quad \quad \quad \forall \gamma \in \Gamma_2(2,4)
\end{equation*}
\ni Hence, $k^2(\gamma)=1$ when $\gamma \in \Gamma$, due to the first identity of the statement; the second identity is therefore proved, $\Gamma^{(-1)}$ being indeed an empty set.\\
\end{proof}

\ni Thanks to Proposition \ref{propcondicio}, a remarkable statement for $\Gamma$ can be proved:

\begin{prop}
$\Gamma$ is normal in $\Gamma_2$.
\end{prop}

\begin{proof}
One has to prove that:
\begin{equation*}
\chi_{N_i} (\gamma^{-1} \eta \gamma) = 1 \quad \quad \forall \gamma \in \Gamma_2 \quad , \quad \forall \eta \in \Gamma \quad , \quad \forall i = 1, \dots , 15
\end{equation*}

\ni An action of $\Gamma_g$ is well defined on the products $\chi_m \chi_n$ by means of the formula $\gamma \cdot \chi_m (\eta) \coloneqq \chi_m (\gamma \eta \gamma^{-1})$ in such a way that $\gamma(\chi_m \chi_n) = (\gamma \cdot \chi_m) (\gamma \cdot \chi_n) = \chi_{\gamma^{-1}m}\chi_{\gamma^{-1}n}$ (cf. \cite{SM}). Therefore, by setting $N_i=(n_{1 \, i},n_{2 \, i})$ for each $i=1, \dots , 15$, one has:
\begin{equation*}
\chi_{N_i} (\gamma^{-1} \eta \gamma) = \chi_{n_{1 \, i}} (\gamma^{-1} \eta \gamma) \chi_{n_{2 \, i}} (\gamma^{-1} \eta \gamma) = \gamma^{-1} (\chi_{n_{1 \, i}}, \chi_{n_{2 \, i}}) (\eta) = \chi_{\gamma n_{1 \, i}} (\eta) \chi_{\gamma n_{2 \, i}} (\eta)
\end{equation*}
Since the action in (\ref{action}) preserves the parity, for each $i=1, \dots , 15$ a $j$ exists, depending on $i$ and $\gamma$, such that $(\gamma n_{1 \, i}, \gamma n_{2 \, i}) = N_j$; hence: 
\begin{equation*}
\chi_{N_j} (\gamma^{-1} \eta \gamma)=\chi_{\gamma n_{1 \, i}} (\eta) \chi_{\gamma n_{2 \, i}}(\eta) = \chi_{N_j} (\eta)=1
\end{equation*}
where the last equality on the right holds by Proposition \ref{propcondicio}, since $\eta \in \Gamma$ as hypothesis.
\end{proof}

\ni A concrete description for $\Gamma$ in terms of generators can be also provided. Since the functions $\chi_{N_i}$ are characters of the group $G$, the elements in $\Gamma_2(2,4)$ belonging to $\Gamma = \bigcap_{i= 1}^{15} \text{Ker}\chi_{N_i}$ can be found by investigating only the representative elements for the cosets of $\Gamma_2(4,8)$ in $\Gamma_2(2,4)$. For such a purpose Proposition \ref{generalgammabasis} will be needed, leading to the following statement in the case $g=2$: 

\begin{lem}
\label{basis}
$G$ is a $9$-dimensional vector space on $\mathds{Z}_2$. A basis is given by:

\begin{equation*}
\begin{array}{ccc} 
        A_{11}=\begin{pmatrix} -1 & 0 & 0 & 0 \\ 0 & 1 & 0 & 0 \\ 0 & 0 & -1 & 0 \\ 0 & 0 & 0 & 1  \end{pmatrix} & A_{12} = \begin{pmatrix} 1 & 2 & 0 & 0 \\ 0 & 1 & 0 & 0 \\ 0 & 0 & 1 & 0 \\ 0 & 0 & -2 & 1  \end{pmatrix} & A_{21} = \begin{pmatrix} 1 & 0 & 0 & 0 \\ 2 & 1 & 0 & 0 \\ 0 & 0 & 1 & -2 \\ 0 & 0 & 0 & 1  \end{pmatrix}\\
\\
        B_{11}^2 = \begin{pmatrix} 1 & 0 & 4 & 0 \\ 0 & 1 & 0 & 0 \\ 0 & 0 & 1 & 0 \\ 0 & 0 & 0 & 1  \end{pmatrix} & B_{22}^2 =\begin{pmatrix} 1 & 0 & 0 & 0 \\ 0 & 1 & 0 & 4 \\ 0 & 0 & 1 & 0 \\ 0 & 0 & 0 & 1  \end{pmatrix} & B_{12} =\begin{pmatrix} 1 & 0 & 0 & 2 \\ 0 & 1 & 2 & 0 \\ 0 & 0 & 1 & 0 \\ 0 & 0 & 0 & 1  \end{pmatrix} \\
\\       C_{11}^2 =\begin{pmatrix} 1 & 0 & 0 & 0 \\ 0 & 1 & 0 & 0 \\ 4 & 0 & 1 & 0 \\ 0 & 0 & 0 & 1  \end{pmatrix} &  C_{22}^2 =\begin{pmatrix} 1 & 0 & 0 & 0 \\ 0 & 1 & 0 & 0 \\ 0 & 0 & 1 & 0 \\ 0 & 4 & 0 & 1  \end{pmatrix} &  C_{12}=\begin{pmatrix} 1 & 0 & 0 & 0 \\ 0 & 1 & 0 & 0 \\ 0 & 2 & 1 & 0 \\ 2 & 0 & 0 & 1  \end{pmatrix}
\end{array}
\end{equation*} 
\end{lem} 

\ni Then, after conventionally enumerating the $6$ odd $2$-characteristics:

\begin{equation*}
\begin{array}{cccccc}
n^{(1)} \coloneqq \begin{bmatrix} 0 1 \\ 0 1 \end{bmatrix} & n^{(2)} \coloneqq \begin{bmatrix} 1 0 \\ 1 0 \end{bmatrix} & n^{(3)} \coloneqq \begin{bmatrix} 0 1 \\ 1 1 \end{bmatrix} &
n^{(4)} \coloneqq \begin{bmatrix} 1 0 \\ 1 1 \end{bmatrix} & n^{(5)} \coloneqq \begin{bmatrix} 1 1 \\ 0 1 \end{bmatrix} & n^{(6)} \coloneqq \begin{bmatrix} 1 1 \\ 1 0 \end{bmatrix}
\end{array} 
\end{equation*} a table can be written down by a straightforward computation, using the values in Lemma \ref{chimonbasis}.

\begin{table}[h]
\label{values}
\centering
\begin{tabular}{c|ccccccccc}
$\chi_{i,j} \coloneqq \chi_{(n^{(i) }, n^{(j)} )}$  & $A_{11}$ & $A_{12}$ & $A_{21}$ & $B_{12}$ & $B^2_{11}$ & $B^2_{22}$ & $C_{12}$ & $C_{11}^2$ & $C_{22}^2$ \\ 
\hline
\\
$\chi_{12}$ & -1 & 1 & 1 & 1 & -1 & -1 & 1 & -1 & -1 \\
$\chi_{13}$ & 1 & 1 & -1 & 1 & 1 & 1 & -1 & -1 & 1 \\
$\chi_{14}$ & -1 & -1 & 1 & 1 & -1 & -1 & -1 & -1 & 1 \\
$\chi_{15}$ & 1 & -1 & 1 & -1 & -1 & 1 & 1 & 1 & 1 \\
$\chi_{16}$ & -1 & 1 & -1 & -1 & -1 & 1 & 1 & -1 & -1 \\
$\chi_{23}$ & -1 & 1 & -1 & 1 & -1 & -1 & -1 & 1 & -1 \\
$\chi_{24}$ & 1 & -1 & 1 & 1 & 1 & 1 & -1 & 1 & -1 \\
$\chi_{25}$ & -1 & -1 & 1 & -1 & 1 & -1 & 1 & -1 & -1 \\
$\chi_{26}$ & 1 & 1 & -1 & -1 & 1 & -1 & 1 & 1 & 1 \\
$\chi_{34}$ & -1 & -1 & -1 & 1 & -1 & -1 & 1 & 1 & 1 \\
$\chi_{35}$ & 1 & -1 & -1 & -1 & -1 & 1 & -1 & -1 & 1 \\
$\chi_{36}$ & -1 & 1 & 1 & -1 & -1 & 1 & -1 & 1 & -1 \\
$\chi_{45}$ & -1 & 1 & 1 & -1 & 1 & -1 & -1 & -1 & 1 \\
$\chi_{46}$ & 1 & -1 & -1 & -1 & 1 & -1 & -1 & 1 & -1 \\
$\chi_{56}$ & -1 & -1 & -1 & 1 & 1 & 1 & 1 & -1 & -1 \\ 
\\
\hline
\end{tabular}
\caption{Values of $\chi_{N_i}$ on a basis of $G$}
\end{table}

\begin{prop}
\label{groupelements}
The group $\Gamma$ is generated by $\Gamma_2(4,8)$ and the elements:
\begin{equation*}
\begin{array}{ccc} 
        A_{12} B_{11}^2 C_{22}^2 & = & \begin{pmatrix} 1 & 2 & 4 & 0 \\ 0 & 1 & 0 & 0 \\ 0 & 0 & 1 & 0 \\ 0 & 4 & -2 & 1  \end{pmatrix} \\
          &   & \\
        ^{t}(A_{12} B_{11}^2 C_{22}^2) =  A_{21} B_{22}^2 C_{11}^2 & = &\begin{pmatrix} 1 & 0 & 0 & 0 \\ 2 & 1 & 0 & 4 \\ 4 & 0 & 1 & -2 \\ 0 & 0 & 0 & 1  \end{pmatrix}\\
          &   & \\
        B_{12} B_{11}^2 B_{22}^2 & = & \begin{pmatrix} 1 & 0 & 4 & 2 \\ 0 & 1 & 2 & 4 \\ 0 & 0 & 1 & 0 \\ 0 & 0 & 0 & 1  \end{pmatrix} \\ 
          &   & \\
        ^{t}(B_{12} B_{11}^2 B_{22}^2) = C_{12} C_{11}^2 C_{22}^2 & = & \begin{pmatrix} 1 & 0 & 0 & 0 \\ 0 & 1 & 0 & 0 \\ 4 & 2 & 1 & 0 \\ 2 & 4 & 0 & 1  \end{pmatrix}
\end{array}
\end{equation*}
\end{prop}

\begin{proof}
Thanks to Proposition \ref{propcondicio} the independent elements satisfying the desired properties can be easily detected by means of the table.
\end{proof}

\begin{coro}
$\mathfrak{S}_2$ does not admit any fixed point under the action of $\Gamma$. In particular, the coset space $\mathfrak{S}_2/\Gamma$ is smooth.
\end{coro}

\begin{proof}
The elements of finite order in $\Gamma_g$ are known to be the only ones fixing points on $\mathfrak{S}_g$; an element in $\Gamma_g(2,4)$ that fixes points on $\mathfrak{S}_g$ has, therefore, order $2$, for $\gamma^2 \in \Gamma_g(2,4)$ whenever $\gamma \in \Gamma_g(4,8)$ and $\Gamma_g(4,8)$ does not fix any point. The thesis then follows, since no element amongst the ones listed in Proposition \ref{groupelements} has order $2$.
\end{proof}

\end{section}

\begin{section}{Structure of $A(\Gamma)$: generators}

\ni The even part $A(\Gamma)^e$ of the ring of modular forms with respect to $\Gamma$ is the only relevant one to describe Satake's compactification $\text{Proj}A(\Gamma)$. A structure theorem have to be proved first, however, in order to describe $A(\Gamma)^e$:
\begin{prop}
\label{a48structure}
Let $\mathds{C}[\theta_m^{2}\theta_n^2]$ be the $\mathds{C}$-algebra generated by the products $\theta_m^{2}\theta_n^2$ with $m$ and $n$ even. Then:
\\
\begin{equation*}
A(\Gamma_2(4,8)) = \left ( \bigoplus_{d=0, 2, 4} \mathds{C}[\theta_m^{2}\theta_n^2] \theta_{m_1} \cdots \theta_{m_{2d}} \right )  \bigoplus \left ( \bigoplus_{d=1, 3, 5} \mathds{C}[\theta_m^{2}\theta_n^2]\theta_{m_1} \cdots \theta_{m_{2d}} \right )
\end{equation*}

\ni where:
\begin{equation*}
A(\Gamma_2(4,8))^e = \bigoplus_{d=0, 2, 4} \mathds{C}[\theta_m^{2}\theta_n^2] \theta_{m_1} \cdots \theta_{m_{2d}}
\end{equation*}
\ni is the even part of the graded ring, and

\begin{equation*}
A(\Gamma_2(4,8))^o = \bigoplus_{d=1, 3, 5} \mathds{C}[\theta_m^{2}\theta_n^2]\theta_{m_1} \cdots \theta_{m_{2d}}
\end{equation*}
is the odd part.

\end{prop}

\begin{proof}
By Igusa's Structure Theorem (cf. \cite{Igusa7} and \cite{Igusa3}), $A(\Gamma_2(4,8) = \mathds{C}[\theta_m \theta_n]$; the ring is thus decomposed under the action of $\Gamma_2(2,4)$ into: 
\begin{equation*}
A(\Gamma_2(4,8))=\bigoplus_{\chi \in \hat{G}}\mathds{C}[\theta_m \theta_n , \chi]
\end{equation*}
\ni where:
\begin{equation*}
\mathds{C}[\theta_m \theta_n , \chi] = \{ P \in \mathds{C}[\theta_m \theta_n] \mid  P(\gamma \tau) = \chi(\gamma) \det{(c_{\gamma} \tau + d_{\gamma})}^k P(\tau) \, \, \forall \gamma \in \Gamma_2(2,4), \, \, k \geq 0 \}
\end{equation*}
\ni Due to (\ref{transformcostanthetagamma2n4n}), one only needs to study how monomials in $\theta_m \theta_n$ transform under the action of $\Gamma_2(2,4)$, to investigate the whole ring. \\ 
\ni One has to note first that whenever $m$ and $n$ are both even or odd characteristics the function $\chi_m \chi_n$ is well defined on $G$, being a character of this group. The set $\{ \chi_m \}_{m \in C_1}$ is, in particular, a set of generators for the group $\hat{G}$ (Cf. \cite{SM}).\\
\ni Let now $P_d=\theta_{m_1}\cdots\theta_{m_{2d}} \in \mathds{C}[\theta_m \theta_n]_d$ be a monomial of degree $d$ in the variables $\theta_m \theta_n$; (\ref{transformcostanthetagamma2n4n}) implies then:
\begin{equation*}
P_d(\gamma \tau) = \kappa^{2d}(\gamma) \chi_{m_1} \cdots \chi_{m_{2d}} \det{(c_{\gamma} \tau + d_{\gamma})}^d P_d(\tau) \quad \quad \forall \gamma \in \Gamma_2(2,4)
\end{equation*}

\ni When $d=2l$, $P_d \in \mathds{C}[\theta_m \theta_n , \chi_{m_1} \cdots \chi_{m_{2d}}]$, for $\kappa^4$ is trivial on $\Gamma(2,4)$. Since each $\chi_m^2\chi_n^2$ is a trivial character of $G$, the following decomposition arises for the even part of the ring:
\begin{equation*}
A(\Gamma_2(4,8))^e = \bigoplus_{d \, \text{even}} \mathds{C}[\theta_m^{2}\theta_n^2] \theta_{m_1} \cdots \theta_{m_{2d}}
\end{equation*} 
\ni When $d=2l+1$, $P_d \in \mathds{C}[\theta_m \theta_n , \kappa^2 \chi_{m_1} \cdots \chi_{m_{2d}}]$. However, $\kappa^2$ is a product of $\chi_m$, the function being an element of $\hat{G}$; hence, $P_d \in \mathds{C}[\theta_m \theta_n ,  \chi_{m_{i_1}} \cdots \chi_{m_{i_r}}  \chi_{m_1} \cdots \chi_{m_{2d}}]$. Therefore:
\begin{equation*}
A(\Gamma_2(4,8))^o = \bigoplus_{d \, \text{odd}} \mathds{C}[\theta_m^{2}\theta_n^2]\theta_{m_1} \cdots \theta_{m_{2d}}
\end{equation*}

\end{proof}

\begin{teo}
\label{agamma}
$A(\Gamma)^e=\mathds{C}[\theta_m^2, D(N)]^{(e)}$. 
\end{teo}

\begin{proof}
By Proposition \ref{propcondicio}, $\Gamma / \{\pm \Gamma_2(4,8)\} \subset G$ is the dual subgroup corresponding to the subgroup $<\chi_{N_1}, \dots , \chi_{N_{15}}> \subset \hat{G}$ generated by the $\chi_{N_i}$; one has, therefore:
\begin{equation*}
A(\Gamma)=\bigoplus_{\chi \in < \chi_{N_i} >}A(\Gamma_2(4,8), \chi)
\end{equation*} where, as above \footnote{Here the symbol $ \mathfrak{O}(\mathfrak{S}_2)$ stands for the space of holomorphic functions on $\mathfrak{S}_2$.}:
\begin{equation*}
A(\Gamma_2(4,8), \chi) = \{ f \in \mathfrak{O}(\mathfrak{S}_2) \,  \mid  f(\gamma \tau) = \chi(\gamma) \det{(c_{\gamma} \tau + d_{\gamma})}^k f(\tau) \, \, \forall \gamma \in \Gamma_2(2,4), \, \, k \geq 0 \}
\end{equation*} Then the thesis follows from Proposition \ref{a48structure}.
\end{proof}

\end{section}

\begin{section}{Structure of $A(\Gamma)$: relations}

The foregoing section has been devoted to describe the generators of $A(\Gamma)^e$. Relations exist amongst these generators, most of which are induced by Riemann's relations; this section aims to provide them through a threefold investigation.

\begin{subsection}{Relations among the $\theta_m^2$}
\label{evenchar}
Relations among $\theta_m^2$ are completely described by Riemann's relations, which can be listed by suitably enumerating the $10$ even $2$-characteristics:
\begin{equation*}
\begin{array}{ccccc}
m^{(1)} \coloneqq \begin{bmatrix} 0 0 \\ 0 0 \end{bmatrix} & m^{(2)} \coloneqq \begin{bmatrix} 0 0 \\ 0 1 \end{bmatrix} & m^{(3)} \coloneqq \begin{bmatrix} 0 0 \\ 1 0 \end{bmatrix} & m^{(4)} \coloneqq \begin{bmatrix} 0 0 \\ 1 1 \end{bmatrix} &  m^{(5)} \coloneqq \begin{bmatrix} 0 1 \\ 0 0 \end{bmatrix} \\
m^{(6)} \coloneqq \begin{bmatrix} 1 0 \\ 0 0 \end{bmatrix} & m^{(7)} \coloneqq \begin{bmatrix} 1 1 \\ 0 0 \end{bmatrix} & m^{(8)} \coloneqq \begin{bmatrix} 0 1 \\ 1 0 \end{bmatrix} & m^{(9)} \coloneqq \begin{bmatrix} 1 0 \\ 0 1 \end{bmatrix} &m^{(10)} \coloneqq \begin{bmatrix} 1 1 \\ 1 1 \end{bmatrix}
\end{array} 
\end{equation*} with $\theta_i \coloneqq \theta_{m^{(i)}}$. There are $15$ biquadratic Riemann's relations:
\begin{equation*}
\begin{array}{ccc}
{\theta}^2_2{\theta}^2_3 = {\theta}^2_1{\theta}^2_4 - {\theta}^2_7{\theta}^2_{10}; &  {\theta}^2_2{\theta}^2_5 = {\theta}^2_7{\theta}^2_9 + {\theta}^2_4{\theta}^2_{8}; & {\theta}^2_3{\theta}^2_5 = {\theta}^2_9{\theta}^2_{10} + {\theta}^2_1{\theta}^2_{8}; \\
 & & \\
{\theta}^2_2{\theta}^2_6 = {\theta}^2_1{\theta}^2_{9} + {\theta}^2_8{\theta}^2_{10}; & {\theta}^2_3{\theta}^2_6 = {\theta}^2_4{\theta}^2_{9} + {\theta}^2_7{\theta}^2_{8}; & {\theta}^2_6{\theta}^2_5 = {\theta}^2_1{\theta}^2_{7} - {\theta}^2_4{\theta}^2_{10}; \\
 & & \\
{\theta}^2_6{\theta}^2_7 = {\theta}^2_3{\theta}^2_{8} - {\theta}^2_1{\theta}^2_{5}; & {\theta}^2_6{\theta}^2_{10} = {\theta}^2_4{\theta}^2_{5} - {\theta}^2_2{\theta}^2_{8}; & {\theta}^2_6{\theta}^2_9 = {\theta}^2_1{\theta}^2_{2} - {\theta}^2_3{\theta}^2_{4}; \\
 & & \\
{\theta}^2_5{\theta}^2_9 = {\theta}^2_2{\theta}^2_{7} - {\theta}^2_3{\theta}^2_{10}; & {\theta}^2_4{\theta}^2_6 = {\theta}^2_5{\theta}^2_{10} + {\theta}^2_3{\theta}^2_{9}; & {\theta}^2_1{\theta}^2_6 = {\theta}^2_5{\theta}^2_{7} - {\theta}^2_2{\theta}^2_{9}; \\
 & & \\
{\theta}^2_6{\theta}^2_8 = {\theta}^2_3{\theta}^2_{7} - {\theta}^2_2{\theta}^2_{10}; & {\theta}^2_5{\theta}^2_8 = {\theta}^2_1{\theta}^2_{3} - {\theta}^2_2{\theta}^2_{4};& {\theta}^2_8{\theta}^2_9 = {\theta}^2_4{\theta}^2_{7} - {\theta}^2_1{\theta}^2_{10};
\end{array}
\end{equation*} which are all independent, and $15$ quartic Riemann's relations, only $5$ of them being independent:
\begin{equation*}
\begin{array}{ccc}
{\theta}^4_1 - {\theta}^4_4 - {\theta}^4_5 - {\theta}^4_9 = 0; & {\theta}^4_2 - {\theta}^4_3 + {\theta}^4_5 - {\theta}^4_6 = 0; & {\theta}^4_2 - {\theta}^4_3 + {\theta}^4_8 - {\theta}^4_9 = 0;\\
 & & \\
{\theta}^4_1 - {\theta}^4_3 - {\theta}^4_6 - {\theta}^4_{10} = 0; & {\theta}^4_1 - {\theta}^4_2 - {\theta}^4_7 - {\theta}^4_8 = 0; & 
\end{array}
\end{equation*}

\ni The $15$ biquadratic Riemann's relations correspond to the elements of $C_6^+$ by the bijective map:
\begin{equation}
\label{rrquad}
M=(m_1, \dots m_6) \quad \longmapsto \quad R_2(M) \quad \, : \, {\theta}^2_{m_1}{\theta}^2_{m_2} \pm {\theta}^2_{m_3}{\theta}^2_{m_4} \pm {\theta}^2_{m_5}{\theta}^2_{m_6}
\end{equation}
\ni The $15$ quartic Riemann's relations correspond to the elements of $C_4^-$ by:
\begin{equation}
\label{rrquart}
M=(m_1, \dots m_4) \quad \longmapsto \quad R_4(M) \quad \, : \, {\theta}^4_{m_1}\pm {\theta}^4_{m_2} \pm {\theta}^4_{m_3} \pm {\theta}^4_{m_4}
\end{equation}\\

\end{subsection}

\begin{subsection}{Relations among the $D(N)$}

\ni By virtue of the general Jacobi's formula (cf. \cite{Igusa1}) the $D(N)$ are monomials of degree $4$ in the Theta constants. More precisely, for each $M=(m_1, m_2, m_3, m_4) \in C_4^{-}$ there exists a Jacobian determinant $D(M)$ such that $D(M) = \theta_{m_1}\theta_{m_2}\theta_{m_3}\theta_{m_4}$, distinct Jacobian determinants being thus set in correspondence with distinct elements of $C_4^{-}$. The bijection $M \mapsto D(M)$ provides, therefore, a parametrization for the Jacobian determinants in the case $g=2$.
\ni With reference to the conventional notations introduced in \S \, \ref{oddchar} and \S \, \ref{evenchar} for odd and even $2$-characteristics, one has, therefore:
\begin{equation*}
\begin{array}{lll}
D(n^{(1)},n^{(2)}) = \theta_{2} \theta_{3} \theta_{5} \theta_{6} ; & D(n^{(1)},n^{(3)}) = - \theta_{6} \theta_{7} \theta_{9} \theta_{10} ; & D(n^{(1)},n^{(4)}) = \theta_{1} \theta_{4} \theta_{5} \theta_{9} \\
 & & \\
D(n^{(1)},n^{(5)}) = - \theta_{3} \theta_{4} \theta_{8} \theta_{10} ; & D(n^{(1)},n^{(6)}) = \theta_{1} \theta_{2} \theta_{7} \theta_{8} ; & D(n^{(2)},n^{(3)}) = - \theta_{1} \theta_{4} \theta_{6} \theta_{8} \\
 & & \\
D(n^{(2)},n^{(4)}) = - \theta_{5} \theta_{7} \theta_{8} \theta_{10}; & D(n^{(2)},n^{(5)}) = - \theta_{1} \theta_{3} \theta_{7} \theta_{9} ; & D(n^{(2)},n^{(6)}) = \theta_{2} \theta_{4} \theta_{9} \theta_{10} \\
 & & \\
D(n^{(3)},n^{(4)}) = \theta_{2} \theta_{3} \theta_{8} \theta_{9} ; & D(n^{(3)},n^{(5)}) = - \theta_{1} \theta_{2} \theta_{5} \theta_{10} ; & D(n^{(3)},n^{(6)}) = \theta_{3} \theta_{4} \theta_{5} \theta_{7} \\
 & & \\
D(n^{(4)},n^{(5)}) = - \theta_{2} \theta_{4} \theta_{6} \theta_{7} ; & D(n^{(4)},n^{(6)}) = \theta_{1} \theta_{3} \theta_{6} \theta_{10} ; & D(n^{(5)},n^{(6)}) = \theta_{5} \theta_{6} \theta_{8} \theta_{9} 
\end{array}
\end{equation*} The relations involving only the $D(N)$ are generated by these ones and by Riemann's relations; here a combinatorial description follows:

\begin{subequations}

\begin{enumerate}
\item For each even characteristic $m$, the six $4$-plets $\{ M_i^m \}_{i=1, \dots 6}$ in $C_4^-$ containing $m$ can be enumerated in such a way that:
\begin{equation*}
M_1^m \cap M_2^m \cap M_3^m = \{ m \} =  M_4^m \cap M_5^m \cap M_6^m 
\end{equation*}
Then, one gains:
\begin{equation}
\label{rb1}
D(M_1^m) D(M_2^m) D(M_3^m) = \chi_5 \theta^2_m  =  D(M_4^m) D(M_5^m) D(M_6^m) 
\end{equation}
which are obviously $10$ relations, namely one for each choice of $m$.

\item For each $M= \{m_1, \dots m_6 \} \in C_6^+$ the eight $4$-plets $\{\tilde{M}_i\}_{i=1, \dots 8}$ in $C_4^-$ containing exactly a triplet $\{m_i, m_j, m_k \} \subset M$ can be enumerated in such a way that:
\begin{equation}
\begin{array}{l}
\label{rb2}
D(\tilde{M}_1)D(\tilde{M}_2) D(\tilde{M}_3)D(\tilde{M}_4) = \chi_5 \prod_{i=1}^6 \theta^2_{m_i}=\\
\\
\quad \quad \quad \quad \quad \quad \quad \quad \quad \quad =D(\tilde{M}_5)D(\tilde{M}_6) D(\tilde{M}_7)D(\tilde{M}_8)
\end{array}
\end{equation}
These are $15$ relations, namely one for each choice of $M \in C_6^+$.

\item Let $M= \{m_1, \dots m_6 \} \in C_6^+$ be fixed and let
\begin{equation*} 
R_2(M) = \theta_{m_1}^2\theta_{m_2}^2 \pm  \theta_{m_3}^2\theta_{m_4}^2 \pm \theta_{m_5}^2\theta_{m_6}^2 = 0
\end{equation*} be the associated biquadratic Riemann's relation as in (\ref{rrquad}). Denoting by $M_1^{i,j}$ and $M_2^{i,j}$ the only two $4$-plets of $C_4^-$ containing $\{m_i, m_j \}$, one has:
\begin{equation}
\begin{array}{l}
\label{rb3}
D(M_1^{1, 2})D(M_2^{1, 2}) \pm D(M_1^{3, 4})D(M_2^{3, 4})  \pm D(M_1^{5, 6})D(M_2^{5,6}) =\\
\\
\quad \quad  \quad \quad  \quad \quad  \quad \quad  \quad \quad \quad \quad \quad = \pm D(M') R_2(M)=0
\end{array}
\end{equation}
\ni These are $15$ relations, corresponding to the elements of $C_6^+$.\\
Otherwise, if only one of the two $4$-plets $M_1^{i,j}$ e $M_2^{i,j}$ is chosen for each couple $\{m_i, m_j \}$, the following general identity is gained:
\begin{equation*}
\begin{array}{l}
D^2(M_{\alpha}^{1, 2}) \pm D^2(M_{\beta}^{3, 4}) \pm D^2(M_{\epsilon}^{5, 6})= \\ 
\\
\quad \quad \quad = {\theta}^2_{\alpha_1}{\theta}^2_{\alpha_2}{\theta}^2_{m_1}{\theta}^2_{m_2} \pm  {\theta}^2_{\beta_1}{\theta}^2_{\beta_2}{\theta}^2_{m_3}{\theta}^2_{m_4} \pm  {\theta}^2_{\epsilon_1}{\theta}^2_{\epsilon_2}{\theta}^2_{m_5}{\theta}^2_{m_6}
\end{array}
\end{equation*}
where $\alpha_1, \alpha_2, \beta_1, \beta_2, \epsilon_1, \epsilon_2$ are characteristics in $M^c \in C_4^-$. However, the triplet of determinants involved in the identity can be chosen in such a way that only three of the four characteristics belonging to $M^c$ appear; for combinatorial reasons, there is a unique way to gain a relation among $D(M)$'s by multiplying each determinant of such a triplet by two other distinct determinants:
\begin{equation}
\label{rb4}
D_h \, D_k \, D^2(M_{\alpha}^{1, 2}) \pm D_l \, D_r \, D^2(M_{\beta}^{3, 4}) \pm D_s \, D_t \, D^2(M_{\epsilon}^{5, 6})=0
\end{equation}
For each $M \in C_6^+$, there exist four triplets $D(M_{\alpha}^{1, 2}), D^2(M_{\beta}^{3, 4}), D^2(M_{\epsilon}^{5, 6})$ satisfying the desired requirement, each triplet corresponding to a choice for the characteristic in $M^c$ which does not appear in the general identity. These relations (\ref{rb4}) are, therefore, $15 \cdot 4 = 60$. 

\item Let $M= \{m_1, \dots m_4 \} \in C_4^-$ and let
\begin{equation*}
R_4(M) = \theta_{m_1}^4 \pm  \theta_{m_2}^4 \pm  \theta_{m_3}^4 \pm \theta_{m_4}^4 = 0
\end{equation*}
be the associated quartic Riemann's relation as in (\ref{rrquart}).
For each $m_i \in M$ there are exactly $2$ elements $M_1^i, M_2^i \in C_4^-$ containing $m_i$ and also satisfying $M_1^i \vartriangle M_2^i = M^c= \{m_5, \dots m_{10} \}$. One has, therefore:
\begin{equation}
\label{rb5}
\sum_{i=1}^4 \pm D(M_1^i)^2 \, D(M_2^i)^2 = \theta_{m_{5}} \cdots \theta_{m_{10}} R_4(M)=0
\end{equation}
The $15$ quartic Riemann's relations all induce independent relations on the $D(M)$, albeit they are not independent themselves; the relations in (\ref{rb5}) are, therefore, $15$.

\item Let $M=\{m_1, m_2, m_3 \} \in C_3^-$ be fixed. \\
Let $\tilde{M}=\{m_1, m_2, m_3, m'_1, m'_2, m'_3 \}$ be one of the two elements in $C_6^+$ containing $M$; the corresponding biquadratic Riemann's relation is:
\begin{equation*}
R_2(\tilde{M}) = \theta_{m_1}^2\theta_{m'_1}^2 \pm  \theta_{m_2}^2\theta_{m'_2}^2 \pm \theta_{m_3}^2\theta_{m'_3}^2 = 0
\end{equation*}
with $\{m'_1, m'_2, m'_3\} \in C_3^-$. Besides, for any couple $\{m_i, m_j \} \subset M$, there exists a unique $M_{i,j} \in C_6^+$, containing $\{m_i, m_j \}$ and satisfying:
\begin{equation*}
R_2(M_{i,j}) = \pm \theta_{m_i}^2\theta_{m_j}^2 + P_{ij} = 0
\end{equation*}
For combinatorial reasons, all the terms $\theta_{m'_i}^2 P_{jk}$ share the common addend $\theta_{m'_1}^2\theta_{m'_2}^2\theta_{m'_3}^2$; therefore, one has:
\begin{equation*}
\begin{split}
0 & = \theta_{m_1}^2\theta_{m_2}^2\theta_{m_3}^2 R_2(\tilde{M}) = \pm  \theta_{m_1}^2\theta_{m'_1}^2 P_{2,3}  \pm  \theta_{m_2}^2\theta_{m'_2}^2 P_{1,3} \pm  \theta_{m_3}^2\theta_{m'_3}^2 P_{1,2} =\\
& \\
& = \pm  \theta_{m_4}^2 \theta_{m'_1}^2\theta_{m'_2}^2\theta_{m'_3}^2  \pm  \theta_{m_1}^2\theta_{m'_1}^2 \theta_{m_{\alpha}}^2 \theta_{m_{\beta}}^2 \pm  \theta_{m_2}^2\theta_{m'_2}^2\theta_{m_{\alpha}}^2 \theta_{m_{\epsilon}}^2  \pm  \theta_{m_3}^2\theta_{m'_3}^2 \theta_{m_{\beta}}^2 \theta_{m_{\epsilon}}^2 
\end{split}
\end{equation*}
where $m_4$ is the unique characteristic completing $M=\{m_1, m_2, m_3 \}$ to an element of $C_4^-$ (to which a quartic Riemann's relations correspond, as in (\ref{rrquart})), and $\{m_{\alpha}, m_{\beta}, m_{\epsilon} \} \coloneqq \{m_1, m_2, m_3, m_4, m'_1, m'_2, m'_3\}^c$. Then:
\begin{equation}
\label{rb6}
0 = \left ( \prod_{m \notin \{m_1, m_2, m_3, m_4 \}} \theta_{m} \right ) \theta_{m_1}^2\theta_{m_2}^2\theta_{m_3}^2 R_2(\tilde{M}) = \sum_{i=1}^4 \pm D(M_1^i) \, D(M_2^i)^3 
\end{equation}
where, for each $i=1, \dots 4$, $M_1^i$ and $M_2^i$ are the $4$-plets in $C_4^-$ containing $m_i$ and such that $M_1^i \vartriangle M_2^i = \{m_1, m_2, m_3, m_4 \}^c$ as in the relations (\ref{rb5}).\\
By choosing the other $6$-plet $\tilde{M}=\{m_1, m_2, m_3, m_{\alpha}, m_{\beta},m_{\epsilon} \}$ containing $M$, the same relation is gained with interchanged exponents:
\begin{equation}
\label{rb7}
0 = \sum_{i=1}^4 \pm D(M_1^i)^3 \, D(M_2^i) 
\end{equation}
Triplets $M= \{m_1, m_2, m_3 \}$ which are contained in the same element of $C_4^-$, yield the same relation (because the same related quartic Riemann's relation turns out to be replaced in the null expression). These relations are thus parameterized by the elements in $C_4^-$; hence, there are $15$ relations of type (\ref{rb6}) and $15$ of type (\ref{rb7}).

\item For each $m \in C_1$, there are exactly six determinants $\{D_i^{m}\}_{i=1, \dots 6}$ such that:
\begin{equation*}
 D_i^{m}=D(M)= \pm \theta_{m_1}\theta_{m_2}\theta_{m_3}\theta_{m_4} \quad \quad \text{with} \, \, m \in M=\{m_1, m_2, m_3, m_4 \} \in C_4^-
\end{equation*}
Hence:
\begin{equation*}
\begin{split}
\sum_{i=1}^6 (D_i^{m})^4  & = \theta^4_m[ \theta^4_{n_1}(\theta^4_{\alpha_1}\theta^4_{\alpha_3} \pm \theta^4_{\alpha_2}\theta^4_{\alpha_4}) \pm \\
& \pm \theta^4_{n_2}(\theta^4_{\alpha_1}\theta^4_{\alpha_5} \pm \theta^4_{\alpha_2}\theta^4_{\alpha_6}) \pm \theta^4_{n_3}(\theta^4_{\alpha_3}\theta^4_{\alpha_5} \pm \theta^4_{\alpha_4}\theta^4_{\alpha_6})]
\end{split}
\end{equation*}
where $\{ \alpha_1, \alpha_2, \alpha_3, \alpha_4\}, \{ \alpha_1, \alpha_2, \alpha_5, \alpha_6\} \in C_4^-$.\\
Now, for each $\{ \alpha_i, \alpha_j, \alpha_k\} \in C_3^-$ there is a unique element $M^k_{ij}$ in $C_4^-$ containing $\alpha_i$, $\alpha_j$, and not $\alpha_k$; by setting $P(M^k_{ij}) \coloneqq R_4(M^k_{ij}) - \theta^4_{\alpha_j}$, where $R_4(M^k_{ij})$ is the Riemann's relation associated to $M^k_{ij}$ as in (\ref{rrquart}), one gains:
\begin{equation*}
\begin{split}
\sum_{i=1}^6 \pm (D_i^{m})^4 & = \theta^4_m \{ \theta^4_{n_1} [\theta^4_{\alpha_1}P(M^1_{23})  \pm \theta^4_{\alpha_2}P(M^2_{14}) ]  \pm \\
& \pm \theta^4_{n_2} [ \theta^4_{\alpha_1}P(M^1_{25}) \pm \theta^4_{\alpha_2}P(M^2_{16}) ] 
\pm \theta^4_{n_3} [\theta^4_{\alpha_3}P(M^3_{45})  \pm \theta^4_{\alpha_4}P(M^4_{36}) ] \}
\end{split}
\end{equation*}
where:
\begin{equation*}
\begin{array}{ccc}
M^1_{23} \cap M^2_{14} = \{n_2, n_3 \}; &  M^1_{25} \cap M^2_{16} = \{n_1, n_3 \}; & M^3_{45} \cap M^4_{36} = \{n_1, n_2 \};
\end{array}
\end{equation*}
By virtue of a suitable choice of the signs, one has then:
\begin{equation*}
\begin{split}
\sum_{i=1}^6 & \pm (D_i^{m})^4 = \theta^4_m [ \theta^4_{n_1}(\theta^4_{\alpha_1} \pm \theta^4_{\alpha_2}) (\pm \theta^4_{n_2} \pm \theta^4_{n_3}) \pm \\
&\theta^4_{n_2}(\theta^4_{\alpha_1} \pm \theta^4_{\alpha_2}) (\pm \theta^4_{n_1} \pm \theta^4_{n_3}) \pm \theta^4_{n_3}(\theta^4_{\alpha_1} \pm \theta^4_{\alpha_2}) (\pm \theta^4_{n_1} \pm \theta^4_{n_2})] = \\
& \\
& =  \theta^4_m (\theta^4_{\alpha_1} \pm \theta^4_{\alpha_2}) [ \theta^4_{n_1} (\pm \theta^4_{n_2} \pm \theta^4_{n_3}) \pm \theta^4_{n_2}(\pm \theta^4_{n_1} \pm \theta^4_{n_3}) \pm \theta^4_{n_3} (\pm \theta^4_{n_1} \pm \theta^4_{n_2})] 
\end{split}
\end{equation*}
which is made null by properly selecting the remaining signs. To sum up, for each $m \in C_1$ one gains the relation:
\begin{equation}
\label{rb8}
\begin{split}
\sum_{i=1}^6 \pm (D_i^{m})^4 = 0
\end{split}
\end{equation}
which is uniquely determined by a suitable choice of the signs. Only $6$ of these ten relations are easily seen to be independent.

\end{enumerate}

\end{subequations}

\begin{prop}
All the relations amongst the $D(N)$ are generated by:

\begin{enumerate} 
\item The $10$ relations in (\ref{rb1});
\item The $15$ relations in (\ref{rb2}); 
\item The $15$ relations in (\ref{rb3}); 
\item The $60$ relations in (\ref{rb4}); 
\item The $15$ relations in (\ref{rb5});
\item The $15$ relations in (\ref{rb6});
\item The $15$ relations in (\ref{rb7});
\item The $6$ relations in (\ref{rb8});
\end{enumerate}

\end{prop}

\begin{proof}
The statement has been proved by means of a Singular program devised by Eberhard Freitag; the relations have been obtained in the ideal generated by the one defined by Riemann's relations and the one defined by the identities $D(N)=\pm \theta_{m_1} \theta_{m_2} \theta_{m_3}\theta_{m_4}$, by using a Gr\"{o}ebner basis to eliminate the $\theta_m$ variables.
\end{proof}

\end{subsection}

\begin{subsection}{Relations between the $D(N)$ and the $\theta_m^2$}

\begin{subequations}
Since the $D(N)$ are modular forms of weight $2$, the relations between the $D(N)$ and the $\theta_m^2$ are indeed between the $D(N)$ and the products $\theta_m^2\theta_n^2$.\\

\ni There are of course the $15$ relations induced by Jacobi's formula:
\begin{equation}
\label{rc1}
D(M)^2=\theta_{m_1}^2\theta_{m_2}^2\theta_{m_3}^2\theta_{m_4}^2 \quad \quad \forall \, M=\{m_1, m_2, m_3, m_4\} \in C_4^-
\end{equation}
\ni Any other relation is clearly generated by (\ref{rc1}) and all the relations of the kind:
\begin{equation*}
\prod_{i = 1}^h D(N_i) = P(\theta_m^2\theta_n^2)
\end{equation*}
\ni where $P(\theta_m^2\theta_n^2)$ is a polynomial in the $\theta_m^2\theta_n^2$ and the $D(N_i)$ are all distinct.

\ni Since for each couple $m,n$ $\theta_m^2\theta_n^2$ is a modular form with respect to $\Gamma_2(2,4)$ (by (\ref{k}) and (\ref{transformcostanthetagamma2n4n})), such a relation holds for a product $\prod_{i = 1}^h D(N_i)$ if and only if $\prod_{i = 1}^h D(N_i)$ is a modular form with respect to $\Gamma_2(2,4)$, this condition being equivalent by Proposition \ref{modular} to $M{}^tM \equiv 0\, \text{mod}\,2$ for the $4 \times 4h$ matrix $M=(M_{1} \dots M_{h})$ where the $M_i$ are the $4$-plets of even characteristics satisfying $D(M_i)=D(N_i)$ (equivalently, to $N{}^tN \equiv 0\, \text{mod}\,2$ for the $4 \times 2h$ matrix $N=(N_{1} \dots N_{h})$ because of (\ref{Nmodular})); a necessary condition is therefore given by $diag(M^tM) \equiv 0 \, \text{mod}\,2$ (or, likewise, by $diag(N^tN) \equiv 0 \,\text{mod}\,2$).

\ni A product $D(N_{1}) \cdots D(N_{h}) = D(M_{1}) \cdots D(M_{h})$ of distinct Jacobian determinants, such that the sum of all the even characteristics $m$ appearing in the $4$plets $M_i \in C_4^-$ (or, equivalently, the sum of all the odd characteristics $n$ appearing in the couples $N_i \in \tilde{C}_2$), each counted with its multiplicity, is congruent to $0$ mod $2$, will be thus called a {\it remarkable factor} of degree $h$. A remarkable factor which is product of remarkable factors will be named {\it reducible}, otherwise it will be called {\it non-reducible}; a product $D(N_{1}) \cdots D(N_{h})$ of distinct determinants which is a modular form with respect to $\Gamma_2(2,4)$ is then a remarkable factor; clearly the converse statement is not necessarily true. Such technical definitions will turn out to be useful, for remarkable factors can be easily characterized:
\begin{prop}
\label{criteriofattori}
$P$ is a remarkable factor if and only if:
\begin{equation*}
P= {\chi_5}^h \prod_{m} {\theta_m}^2 \quad \quad h=0, 1
\end{equation*} where $\chi_5 \coloneqq \chi^{(2)} \coloneqq \prod_{m \in C_1} \theta_m$ is the product of the $10$ non trivial Theta constants.
\end{prop}

\begin{proof}
A product $P$ of distinct $D(N)$, which admits an expression as a monomial in $\theta^2_m$ and $\chi_5$, is clearly a remarkable factor; to prove the converse statement one can define the following function on the set $P(C^{(2)})$ of the parts of $C^{(2)}$: \\
\begin{equation*}
\begin{array}{ccc}
\begin{array}{l}
F: P(C^{(2)})  \longrightarrow \mathds{C}[\theta_m]\\
\\
 \quad  \{m_1, m_2, \dots , m_h \}  \longrightarrow   \theta_{m_1}\theta_{m_2} \cdots \theta_{m_h}
\end{array} & \quad \quad & 
F (\emptyset) \coloneqq 1
\end{array}
\end{equation*}
 
\ni Clearly $F(\{m\})=\theta_m$, $ F(C_1)=\chi_5$, and $F(M)=D(M)$ when $M \in C_4^{-}$. Then, since:
\begin{equation*}
F(M_i)\,F(M_j)=  F(M_1 \vartriangle M_j) \prod_{m \in M_i \cap M_j} {\theta_m}^2
\end{equation*}
\ni if $P=F(M_1) \cdots F(M_h)$ is a remarkable factor, Propositions \ref{pA}, \ref{pB} and \ref{pC} imply $P= \prod_{m} {\theta_m}^2$ or $P= \chi_5 \prod_{m} {\theta_m}^2$.
\end{proof}

\ni Remarkable factors can be classified by means of the law:
\begin{equation*}
N=\{n_1, n_2 \} \, \, \longrightarrow \, \, S(N) \coloneqq n_1 + n_2
\end{equation*} associating to each couple of odd characteristics $N \in \tilde{C}_2$ their sum $S(N) \in \mathds{Z}_2^4$. A product $P=\prod_i D(N_i)$ of distinct Jacobian determinants is, in particular, a remarkable factor if and only if $\sum_iS(N_i)=0$.

\begin{lem}
\label{morethan5}
Remarkable factors of degree greater than $5$ are reducible. 
\end{lem}

\begin{proof}
Let $P= \prod_{i=1}^hD(N_i)$ with $h >5$ be a remarkable factor. \\
The set $\{S(N_i)\}_{i=1, \dots h} \subset \mathds{Z}_2^4$ contains at least two elements linearly dependent from the others. Since $S(N) \neq 0$ for each $N \in \tilde{C}_2$, the thesis follows.
\end{proof}

\ni By Lemma \ref{morethan5}, the remarkable factors whose degree is at most $5$ are the only ones to be checked, in order to find all the non-reducible remarkable factors:

\begin{prop}
\label{listnonred}
The non-reducible remarkable factors are:
\begin{enumerate}
\item $D(n_i, n_j) D(n_j, n_k) D(n_k, n_i)$; 
\item $D(n_i, n_j) D(n_k, n_l) D(n_s, n_t)$;
\item $D(n_i, n_j) D(n_j, n_k) D(n_k, n_l) D(n_l, n_i)$;
\item $D(n_i, n_j) D(n_i, n_k) D(n_i, n_l) D(n_s, n_t)$;
\item $D(n_i, n_j) D(n_j, n_k) D(n_k, n_l) D(n_l, n_r) D(n_r, n_i)$;
\item $D(n, n_j) D(n, n_k) D(n, n_l) D(m, n_r) D(m, n_s)$;
\end{enumerate}
\end{prop}

\begin{proof}
They can be plainly detected by means of Table 2:

\begin{table}[h]
\label{SN}
\centering
\begin{tabular}{|c|c|}
\hline
$D_{ij} \coloneqq D(n^{(i)}, n^{(n_j)})$ & ${}^tS(N)$ \\ \hline
$D_{12}$ & $(1 1 1 1)$ \\
$D_{13}$ & $(0 0 1 0)$ \\
$D_{23}$ & $(1 1 0 1)$ \\
$D_{14}$& $(1 1 1 0)$ \\
$D_{24}$ & $(0 0 0 1)$ \\
$D_{15}$ & $(1 0 0 0)$ \\
$D_{25}$ & $(0 1 1 1)$ \\
$D_{16}$& $(1 0 1 1)$ \\
$D_{26}$ & $(0 1 0 0)$ \\
$D_{34}$& $(1 1 0 0)$ \\
$D_{35}$ & $(1 0 1 0)$ \\
$D_{36}$ & $(1 0 0 1)$ \\
$D_{45}$& $(0 1 1 0)$ \\
$D_{46}$ & $(0 1 0 1)$ \\
$D_{56}$ & $(0 0 1 1)$ \\ 
\hline
\end{tabular}
\caption{Values of S(N)}
\end{table}

\end{proof}

\ni Only factors of type $2$, $3$ and $6$ are modular forms with respect to $\Gamma_2(2,4)$; by Proposition \ref{criteriofattori}, no Theta constant appears in such factors with even multiplicity, while a $\chi_5$ appears with odd multiplicity in factors of type $1$, $4$ and $5$.\\ 

All The products $\prod_iD(N_i)$ of distinct determinants which are monomials in $\theta_m^2\theta_n^2$, are, therefore, products of the factors listed in Proposition \ref{listnonred}.

\begin{prop}
\label{3products}
The relations involving products of $3$ determinants are:
\begin{equation}
\label{rc2}
D(n_i, n_j) D(n_k, n_l) D(n_s, n_t) = \prod_{i=1}^6{\theta_{m_i}}^2
\end{equation}
\end{prop}

\begin{proof}
Factors of type $2$ are the only ones involved. To prove the six Theta constants appearing in the expression (\ref{rc2}) are all distinct, let $P$ be a product of determinants as in (\ref{rc2}) and $M_1, M_2, M_3 \in C_4^-$ the $4$-plets satisfying $P = D(M_1)D(M_2)D(M_3)$. Since $P$ is a remarkable factor, if $M_i \vartriangle M_j \in C_6^+$ holded for any couple of this $4$-plets, it would imply $(M_1 \vartriangle M_2)^c =M_3$; then, one would have $P={\theta_{m_1}}^2 {\theta_{m_2}}^2 {\theta_{m_3}}^2 \chi_5$, which is an absurd statement, for $P$ is also a modular form with respect to $\Gamma_2 (2,4)$. Hence, by Proposition \ref{pA}, $M_i \vartriangle M_j \in C_4^{-}$ for each distinct couple $M_i, M_j$, and the only possibility is $M_1 \vartriangle M_2 = M_3 \in C_4^{-}$, namely: 
\begin{equation*}
M_1= \{ m_1, m_2, m_3, m_4 \} \quad \quad \quad M_2= \{ m_1, m_2, m_5, m_6 \} \quad \quad \quad  M_3 = \{ m_3, m_4, m_5, m_6 \}
\end{equation*} which concludes the proof.
\end{proof}

\begin{prop}
\label{4products}
The relations involving products of $4$ determinants are:
\begin{equation}
\label{rc3}
D(n_i, n_j) D(n_j, n_k) D(n_k, n_l) D(n_l, n_i) =  \prod_{i=1}^8{\theta_{m_i}}^2
\end{equation}

\end{prop}

\begin{proof}
Factors of type $3$ are the only ones involved. To prove the eight Theta constants appearing in the expression (\ref{rc3}) are all distinct, let $P$ be a product of determinants as in (\ref{rc3}) and $M_1, M_2, M_3, M_4 \in C_4^-$ the $4$-plets satisfying $P=D(M_1)D(M_2)D(M_3)D(M_4)$.\\
If $M_1 \vartriangle M_2 \in C_4^-$ , then $M_3 \vartriangle M_4= M_1 \vartriangle M_2 \in C_4^{-}$, for $P$ is a modular form with respect to $\Gamma_2 (2,4)$. Since $|M_1 \vartriangle  M_2| = 4$, there are at least six distinct characteristics appearing in the expression, each with multiplicity $2$; however, $M_3 \vartriangle M_4 = M_1 \vartriangle M_2$ and $|M_3 \cap M_4| = 2$, hence the eight characteristics appearing with multiplicity $2$ are all distinct, for the determinants involved $D(M_i)$ are distinct.\\
If $M_1 \vartriangle M_2 \in C_6^+$, then $M_3 \vartriangle M_4= M_1 \vartriangle M_2 \in C_6^{+}$. Since $|M_1 \cap M_2| = 1$, at least seven distinct characteristics appear, each with multiplicity $2$.  As before, $M_3 \vartriangle M_4= M_1 \vartriangle M_2$ with $|M_3 \cap M_4| = 1$ and the common characteristic in $M_3$ and $M_4$ must be different from the other seven, since the $D(M_i)$ are all distinct.
\end{proof}

\begin{prop}
\label{5products}
The relations involving products of $5$ determinants are:
\begin{equation}
\label{rc4}
D(n, n_j) D(n, n_k) D(n, n_l) D(m, n_r) D(m, n_s)= \prod_{m}\theta_m^2
\end{equation}
\end{prop}

\begin{proof}
Factors of type $6$ are the only ones involved. By virtue of Proposition \ref{criteriofattori}, $\chi_5$ appears with even multiplicity, but the ten Theta constants in the expression (\ref{rc4}) are plainly seen to be not necessarily all distinct.
\end{proof}

\ni As concerns the relations involving products of more than $5$ determinants, the product of two non-reducible remarkable factors of type $2$, $4$ and $5$ is indeed a modular form with respect to $\Gamma_2(2,4)$; therefore, if such a product does not split into some of the factors which have been listed above, it will induce new independent relations. One has, in particular, the following:

\begin{prop}
\label{6products}
The relations involving products of $6$ determinants are:
\begin{equation}
\label{rc5}
D(n_i, n_j) D(n_j, n_k) D(n_k, n_i) D(n_l, n_r) D(n_r, n_s) D(n_s, n_l)= {\chi_5}^2 {\theta_{m}}^2 {\theta_{n}}^2
\end{equation}

\end{prop}

\begin{proof}
The only possible case rises from the product of two distinct factors of type $1$:
\begin{equation*}
\begin{array}{l}
Q_1 = D(n_i, n_j) D(n_j, n_k) D(n_k, n_i) = {\chi_5} \theta_{m}^2\\
\\
Q'_1 = D(n'_i, n'_j) D(n'_j, n'_k) D(n'_k, n'_i) = {\chi_5}\theta_{n}^2
\end{array}
\end{equation*}
Since $Q_1 \cdot Q'_1$ does not factorize into products of determinants which are in turn modular forms with respect to $\Gamma_2(2,4)$, the relations in (\ref{rc5}) are independent from the ones listed above.
\end{proof}

\end{subequations}

\ni The following Proposition ends the investigation around these relations.

\begin{prop}
Let $P$ be a product of more than $6$ distinct determinants, which is a modular form with respect to $\Gamma_2(2,4)$. Then, each relation involving $P$ is dependent from the ones in (\ref{rc1}), (\ref{rc2}), (\ref{rc3}), (\ref{rc4}) and (\ref{rc5}).
\end{prop}

\begin{proof}
The single cases are to be briefly discussed.\\
Let $P$ be a product of $7$ distinct determinants such that $P \in A(\Gamma_2(2,4))$.Then $P$ is necessarily the product of a factor $P_1$ of type $1.$ and a factor $P_4$ of type $4.\,$ and the only possible cases are:
\begin{equation*}
\begin{array}{l}
P_1 \cdot P_4 = [D(n_i, n_j) D(n_j, n_k) D(n_k, n_i)] [D(n_l, n_i) D(n_l, n_j) D(n_l, n_k) D(n_r, n_s)] \\
\\
P_1 \cdot P_4 = [D(n_i, n_j) D(n_j, n_k) D(n_k, n_i)] [D(n_l, n_i) D(n_l, n_j) D(n_l, n_r) D(n_k, n_s)]
\end{array}
\end{equation*}

\ni However, by using the relations (\ref{rb1}), it turns out that:
\begin{equation*}
 D(n_i, n_j) D(n_j, n_k) D(n_k, n_i) =  D(n_l, n_r) D(n_r, n_s) D(n_s, n_l) \\
\end{equation*}

\ni In both cases at least a $D(N)^2$ appears; hence, the relations involving $P$ are dependent from the ones already found (the relations in (\ref{rc1}) hold, in particular).\\
Concerning with products $P=\prod_{N \in C \subset \tilde{C}_2} D(N)$ involving more than $7$ determinants, the product of the determinants which are associated to the complementary couples, namely $P^c \coloneqq \prod_{N \notin C}D(N)$, can be more easily investigated. In fact, when $P \in A(\Gamma_2(2,4))$, then $P^c \in A(\Gamma_2(2,4))$; hence $P^c$ pertains to the previous cases.\\
Let $P$ be, then, a product of $8$ distinct determinants such that $P \in A(\Gamma_2(2,4))$. Since $P^c$ has degree $7$, either $P^c=Q_1 \cdot Q_4$ with $Q_1$ of type $1$ and $Q_4$ of type $4$, or $P^c=Q_2 \cdot Q_3$ with $Q_2$ of type $2$ and $Q_3$ of type $3$.

\ni If $P^c = Q_1 \cdot Q_4$ the only possible cases have been previously discussed; then, it is easily verified that $P$ always admits a factor of the type (\ref{rc4}), which does not appear in the product $Q_1 \cdot Q_4$. Therefore, $P$ splits into the product of two factors which are modular forms with respect to $A(\Gamma_2(2,4))$, these ones having been previously investigated.

\ni If $P^c=Q_2 \cdot Q_3$, then: 
\begin{equation*}
P^c = D(n_i, n_j) D(n_k, n_l) D(n_r, n_s) D(n'_i, n'_j) D(n'_j, n'_k) D(n'_k, n'_l) D(n'_l, n'_i) \\
\end{equation*}
Since four of the six odd characteristics appear with multiplicity $3$ and the other two with multiplicity $1$, $P$ always contains a factor of the type (\ref{rc3}); therefore, $P$ factorizes anyway into the product of two factors which are modular forms with respect to $A(\Gamma_2(2,4))$.\\
Now, let $P$ be the product of $9$ distinct determinants such that $P \in A(\Gamma_2(2,4))$. Since $P^c \in A(\Gamma_2(2,4))$ has degree $6$, it is either the product $Q_2 \cdot Q'_2$ of two factors of type $2$, or the product $Q_1 \cdot Q_1'$ of two factors of type $1$. In the first case:
\begin{equation*}
P^c = Q_2 \cdot Q'_2 =  D(n_i, n_j) D(n_k, n_l) D(n_r, n_s) D(n'_i, n'_j) D(n'_k, n'_l) D(n'_r, n'_s)\\
\end{equation*}
Hence, $P$ always contains a factor of type (\ref{rc4}), all the characteristics being involved with multiplicity $2$.In the second case $P^c$ is of the type (\ref{rc5}); then, at least five characteristics appear in $P^c$ with multiplicity $2$. Then, $P$ always contains a factor of type (\ref{rc4}).\\
\ni If $P$ is a product of $10$ distinct determinants such that $P \in A(\Gamma_2(2,4))$, $P^c$ is of type (\ref{rc4}), by Proposition \ref{5products}. Then, $P$ is easily seen to contain a factor of type (\ref{rc3}).\\
If $P$ is a product of  $11$ distinct determinants such that $P \in A(\Gamma_2(2,4))$, $P^c$  is of type (\ref{rc3}) by Proposition \ref{4products}. Then, $P=P_ 6 \cdot P' $ where $P_6$ is of type (\ref{rc4}) and $P'=P_1 \cdot P_1$ is of type (\ref{rc5}).\\
Finally, if $P$ is a product of $12$ distinct determinants such that $P \in A(\Gamma_2(2,4))$, $P^c$ is of type (\ref{rc2}) by Proposition \ref{3products}. Then, $P=P_ 6 \cdot P_3 \cdot P_2 $ where $P_6$ is of type (\ref{rc4}), $P_3$ is of type (\ref{rc3}) and $P_2$ is of type (\ref{rc2}).\\

\ni The statement is thus proved, since products of $13$ or $14$ Jacobian determinants can not be modular forms with respect to $\Gamma_2(2,4)$, while the product of all the $15$ determinants trivially factorizes into factors which have been already investigated. 

\end{proof}

\ni To sum up, the following Proposition can be stated:
\begin{prop}
A system of independent relations between $D(N)$ and $\theta_m^2\theta_n^2$ is given by (\ref{rc1}), (\ref{rc2}), (\ref{rc3}), (\ref{rc4}) and (\ref{rc5}). 
\end{prop}

\end{subsection}

A description of $A(\Gamma)^{(e)}$ in terms of relations is thus provided by gathering all the results proved in the section:
\begin{teo}
The ideal of the relations amongst $D(N)$ and $\theta_m^2$ is generated by the $20$ independent Riemann's relations, plus the relations (\ref{rb1}), (\ref{rb2}), (\ref{rb3}), (\ref{rb4}), (\ref{rb5}),  (\ref{rb6}), (\ref{rb7}), (\ref{rb8}), (\ref{rc1}), (\ref{rc2}), (\ref{rc3}), (\ref{rc4}) and (\ref{rc5}).
\end{teo}

\end{section}

\begin{section}{The Ideal $S(\Gamma)$}
Since the $D(N)$ are cusp forms, the map $\mathds{P}grTh_2$ does not extend to the boundary of Satake's compactification. In order to describe the desingularization $\text{Proj}S(\Gamma)$, this ending section aims to describe the even part $S(\Gamma)^{e}$ of the ideal of cusp forms with respect to the subgroup $\Gamma$. More precisely, the following structure theorem holds:

\begin{teo}
A system of generators for $S(\Gamma)^e$ is given by:
\begin{equation*}
\begin{array}{ll}
1.& D(M) \quad \quad \forall \, M \in C_4^- \\
    & \\
2. &  \theta_{m_1}^4 \theta_{m_2}^2 \theta_{m_3}^2 \theta_{m_4}^2 \theta_{m_5}^2  \quad \quad \forall \, \{m_1, m_2, m_3, m_4, m_5 \} \in C_5^*
\end{array}
\end{equation*}
In particular, there are $15 + 5 \cdot 72 = 375$ generators for this ideal.
\end{teo}

\begin{proof}
Since $S(\Gamma) \subset S(\Gamma_2(2,4,8))$, the generators of $S(\Gamma)^e$ are amongst the ones described in \cite{Structure} (Theorem 4.4). $S(\Gamma)^e \subset \mathds{C}[\theta_m^2\theta_n^2, D(M)]$ by Theorem \ref{agamma}, hence only types $1.$ and $2.$ in the statement generates $S(\Gamma)^e$; in fact, by using the relations $\theta^2_m=Q_m(\Theta_{m'})$, involving Theta constants and second order Theta constants, namely $\Theta_{m'}(\tau) \coloneqq \theta_{\left [\substack{m' \\ 0} \right ]}(2\tau)$ (cf. \cite{Structure}), the $240$ further elements generating $S(\Gamma_2(2,4,8))$ are easily seen no to belong to $S(\Gamma)^e$, since they can be expressed as $P(\theta_m^2\theta_n^2)\,\Theta_{m'}$, where $P(\theta_m^2\theta_n^2)$ is a polynomial in $\theta_m^2\theta_n^2$.
\end{proof}

\end{section}


\addcontentsline{toc}{chapter}{Bibliography}

\begin{thebibliography}{9}
\bibitem[CS03a]{capsern1}{\it L. Caporaso, E. Sernesi}, Recovering plane curves form their bitangents, Journal of Algebraic Geometry n.12 (2003)
\bibitem[CS03b]{capsern2}{\it L. Caporaso, E. Sernesi}, Characterizing curves by their odd theta-characteristics, arXiv:math/0204164v2 (2003)
\bibitem[GSM03]{mainarticle} {\it S. Grushevsky, R. Salvati Manni}, Gradients of odd theta functions (2003) 
\bibitem[GS93]{Structure}{\it B. Van Geemen, D. Van Straten}, The cusp forms of weight $3$ on $\Gamma_2(2,4,8)$, Mathematics of Computation 61 (1993) 
\bibitem[I64a]{Igusa7} {\it J. Igusa}, On Siegel modular forms of genus $2$ (II), American Journal of Mathematics 86 (1964)
\bibitem[I64b]{Igusa3} {\it J. Igusa}, On the graded ring of theta-constants, American Journal of Mathematics 86 (1964)
\bibitem[I66]{Igusa10}{\it J. Igusa}, On the graded ring of theta-constants (II), American Journal of Mathematics 88 (1966)
\bibitem[I80]{Igusa1} {\it J. Igusa}, On Jacobi's derivative formula and its generalizations, American Journal of Mathematics 102 (1980) n.2
\bibitem[SM83]{SalvatiManniPhd} {\it R. Salvati Manni}, On the nonidentically zero Nullwerte of Jacobians of theta functions with odd characteristics, Adv. in Math. 47 (1983) n.1
\bibitem[SM94]{SM} {\it R. Salvati Manni}, Modular varieties with level 2 Theta strucutre, American Journal of Mathematics 116 (1994) n.6
\bibitem[SMT93]{CuspForms} {\it R. Salvati Manni, J.Top}, Cusp Forms of weight 2 for the group $\Gamma_2(4,8)$, American Journal of Mathematics (1993)
\end{thebibliography}
\end{document}